\documentclass[12pt,reqno]{amsart}
\usepackage[czech, english]{babel}
\usepackage{url}

\newcommand{\tf}{\tfrac}

 \renewcommand{\a}{\alpha}
\renewcommand{\b}{\beta}

\newcommand{\G}{\Gamma}

\renewcommand{\(}{\left\(}
\renewcommand{\)}{\right\)}
\renewcommand{\[}{\left\[}
\renewcommand{\]}{\right\]}

\numberwithin{equation}{section}
 \theoremstyle{plain}
\newtheorem{theorem}{Theorem}[section]
\newtheorem{lemma}[theorem]{Lemma}
\newtheorem{corollary}[theorem]{Corollary}

   \makeatletter
\def\proof{\@ifnextchar[{\@oproof}{\@nproof}}
\def\@oproof[#1][#2]{\trivlist\item[\hskip\labelsep\textit{#2 Proof of\
#1.}~]\ignorespaces}
\def\@nproof{\trivlist\item[\hskip\labelsep\textit{Proof.}~]\ignorespaces}
\makeatother

\begin{document}
\title[An integral analogue of partial theta function]{Error functions, Mordell integrals and an integral analogue of partial theta function}
\author{Atul Dixit}
\address{Department of Mathematics, Tulane University, New Orleans, LA 70118, USA}
\email{adixit@tulane.edu}
\author{Arindam Roy}
\thanks{2010 \textit{Mathematics Subject Classification.} Primary 11M06, 33B20; Secondary 33C15, 34E05.\\
\textit{Keywords and phrases.} Mordell integrals, Ramanujan's Lost Notebook, partial theta function, Riemann zeta function, Dirichlet $L$-function, hypergeometric functions, asymptotic expansions.
}
\address{Department of Mathematics, University of Illinois, 1409 West Green
Street, Urbana, IL 61801, USA} \email{roy22@illinois.edu}
\author{Alexandru Zaharescu}
\address{Department of Mathematics, University of Illinois, 1409 West Green
Street, Urbana, IL 61801, USA and Simion Stoilow Institute of Mathematics of the
Romanian Academy, P.O. Box 1-764, RO-014700 Bucharest, Romania} \email{zaharesc@illinois.edu}
\dedicatory{\emph{Dedicated to Professor Bruce C.~Berndt on the occasion of his 75th birthday}}
\begin{abstract}
A new transformation involving the error function $\textup{erf}(z)$, the imaginary error function $\textup{erfi}(z)$, and an integral analogue of a partial theta function is given along with its character analogues. Another complementary error function transformation is also obtained which when combined with the first explains a transformation in Ramanujan's Lost Notebook termed by Berndt and Xu as the one for an integral analogue of theta functions. These transformations are used to obtain a variety of exact and approximate evaluations of some non-elementary integrals involving hypergeometric functions. Several asymptotic expansions, including the one for a non-elementary integral involving a product of the Riemann $\Xi$-function of two different arguments, are obtained, which generalize known results due to Berndt and Evans, and Oloa.
\end{abstract}
\maketitle
\section{Introduction}
Mordell initiated the study of the integral $\displaystyle\int_{-\infty}^{\infty}\frac{e^{ax^2+bx}}{e^{cx}+d}\, dx$, Re$(a)<0$, in his two influential papers \cite{mordell1, mordell2}. Prior to his work, special cases of this integral had been studied, for example, by Riemann in his Nachlass \cite{siegel} to derive the approximate functional equation for the Riemann zeta function, by Kronecker \cite{kronecker1, kronecker2} to derive the reciprocity for Gauss sums, and by Lerch \cite{lerch1, lerch2, lerch3}. Mordell showed that the above integral can be reduced to two standard forms, namely,
\begin{align}
\varphi(z,\tau)&:=\tau\int_{-\infty}^{\infty}\frac{e^{\pi i\tau x^2-2\pi zx}}{e^{2\pi x}-1}\, dx,\label{mint1}\\
\sigma(z,\tau)&:=\int_{-\infty}^{\infty}\frac{e^{\pi i\tau x^2-2\pi zx}}{e^{2\pi \tau x}-1}\, dx,\nonumber
\end{align}
for Im$(\tau)>0$, and was the first to study the properties of these integrals with respect to modular transformations. Bellmann \cite[p.~52]{bellman} coined the terminology `Mordell integrals' for these types of integrals.

Mordell integrals play a very important role in the groundbreaking Ph.D. thesis of Zwegers \cite{zwegers} which sheds a clear light on Ramanujan's mock theta functions. The definition of a Mordell integral $h(z;\tau)$, Im$(\tau)>0$, employed by Zwegers \cite[p.~6]{zwegers}, and now standard in the contemporary literature is,
\begin{equation*}
h(z;\tau):=\int_{-\infty}^{\infty}\frac{e^{\pi i\tau x^2-2\pi zx}}{\cosh\pi x}\, dx.
\end{equation*}
As remarked by Zwegers himself \cite[p.~5]{zwegers}, $h(z;\tau)$ is essentially the function $\varphi(z;\tau)$ defined in \eqref{mint1}, i.e.,
\begin{equation}\label{relvph}
h(z;\tau)=\tfrac{-2i}{\tau}e^{-\left(\frac{\pi i\tau}{4}+\pi iz\right)}\varphi\left(z+\tfrac{\tau-1}{2},\tau\right).
\end{equation}
Kuznetsov \cite{kuznetsov} has recently used $h(z;\tau)$ to simplify Hiary's algorithm \cite{hiary2} for computing the truncated theta function $\sum_{k=0}^{n}\text{exp}(2\pi i(zk+\tau k^2))$, which in turn is used to compute $\zeta\left(\frac{1}{2}+it\right)$ to within $\pm t^{-\lambda}$ in $O_{\lambda}(t^{\frac{1}{3}}(\ln t)^{\kappa})$ arithmetic operations \cite{hiary1}. We refer the reader to a more recent article \cite{cr} and the references therein for further applications of Mordell integrals.

In \cite{ram1915} and \cite{ram1919}, Ramanujan studied the integrals
\begin{align*}
\phi_{\omega}(z):=\int_{0}^{\infty}\frac{\cos(\pi xz)}{\cosh(\pi x)}e^{-\pi\omega x^2}\, dx,\\
\psi_{\omega}(z):=\int_{0}^{\infty}\frac{\sin(\pi xz)}{\sinh(\pi x)}e^{-\pi\omega x^2}\, dx.
\end{align*}
Of course, we require Re $\omega>0$ for the integrals to converge. If we replace $\omega$ by $-i\tau$ with Im$(\tau)>0$ and $z$ by $2iz$, then the integral $\phi$ is nothing but the Mordell integral. That is,
\begin{equation*}
h(z,\tau)=2\phi_{-i\tau}(2iz).
\end{equation*}
Later, Ramanujan briefly worked on these two integrals in a two-page fragment transcribed by G.~N.~Watson from Ramanujan's loose papers and published along with Ramanujan's Lost Notebook \cite[p.~221-222]{lnb}. See also \cite[pp.~307-328]{geabcbrln4} for details. 

A third integral of this kind studied by Ramanujan on page 198 of the Lost Notebook is
\begin{equation*}
F_{\omega}(z):=\int_{0}^{\infty}\frac{\sin(\pi xz)}{\tanh{(\pi x)}}e^{-\pi \omega x^2}\, dx.
\end{equation*}
As before, one needs Re $\omega>0$ for convergence. One can easily rephrase this integral as
\begin{equation}\label{ftzrep}
F_{\omega}(z)=\int_{-\infty}^{\infty}\frac{e^{-\pi\omega x^2}\sin(\pi xz)}{e^{2\pi x}-1}\, dx=\int_{-\infty}^{\infty}\frac{e^{-\pi\omega x^2+2\pi x}\sin(\pi xz)}{e^{2\pi x}-1}\, dx.
\end{equation}
For Im$(\tau)>0$, this integral $F$ is connected to the integral $\varphi(z,\tau)$ in \eqref{mint1} (and hence to the Mordell integral $h(z;\tau)$) via
\begin{equation*}
F_{-i\tau}(2iz)=\frac{1}{2i\tau}\left(\varphi(z,\tau)-\varphi(-z,\tau)\right).
\end{equation*}
Thus Mordell integrals pervade Ramanujan's papers and his Lost Notebook. In a further support to this claim, we refer the readers to an interesting paper of Andrews \cite{andmor}.

Berndt and Xu \cite{berndtxu} have proved all of the properties of $F_{\omega}(z)$ claimed by Ramanujan in the Lost Notebook. Following Ramanujan, we assume $\omega>0$. Suppose a certain result holds for $\omega>0$, it is clear that by analytic continuation, one may be able to extend it to complex values of $\omega$ in a certain region containing the positive real line. Among the various properties claimed by Ramanujan is the transformation
\begin{equation}\label{fwttrans}
F_{\omega}(z)=\frac{-i}{\sqrt{\omega}}e^{-\frac{\pi z^2}{4\omega}}F_{1/\omega}\left(\frac{iz}{\omega}\right).
\end{equation}
In view of this property, Berndt and Xu call $F_{\omega}(z)$ as an integral analogue of a theta function.

Assume $\alpha>0$, let $\omega=\alpha^2$ and replace $z$ by $\a z/\sqrt{\pi}$ in \eqref{fwttrans}. Using \eqref{ftzrep}, one sees that the above transformation translates into the following identity.

\textit{For $\alpha, \beta>0$ such that $\a\b=1$, 
\begin{align}\label{rt}
\sqrt{\alpha}e^{\frac{z^2}{8}}\int_{-\infty}^{\infty}\frac{e^{-\pi\alpha^2 x^2}\sin(\sqrt{\pi}\alpha xz)}{e^{2\pi x}-1}\, dx=\sqrt{\beta}e^{\frac{-z^2}{8}}\int_{-\infty}^{\infty}\frac{e^{-\pi\beta^2 x^2}\sinh(\sqrt{\pi}\beta xz)}{e^{2\pi x}-1}\, dx.
\end{align}}

Now consider the integral 
\begin{equation}\label{intanpar}
\int_{0}^{\infty}\frac{e^{-\pi\alpha^2 x^2}\sin(\sqrt{\pi}\alpha xz)}{e^{2\pi x}-1}\, dx.
\end{equation}
In analogy with the partial theta function which is the same as a theta function but summed over only half of the integer lattice, we call the above integral an integral analogue of partial theta function. Note that unlike $F_{\omega}(z)$, this integral cannot be expressed solely in terms of Mordell integrals. Also unlike how the transformation formula for the Jacobi theta function trivially gives the transformation formula for the corresponding partial theta function \cite[p.~22, Equation (2.6.3)]{titch}, the transformation in \eqref{rt} does not give rise to a corresponding transformation for the integral in \eqref{intanpar}. (It is easy to check that the integrands in \eqref{rt} are not even functions of $x$). Nevertheless, the primary goal of this paper is to prove a new and interesting modular transformation for the integral in \eqref{intanpar} which involves error functions. 

The error function $\textup{erf}(z)$ and the complementary error function $\textup{erfc}(z)$, defined by \cite[p.~275]{temme}
\begin{equation}\label{erfz}
\textup{erf}(z)=\frac{2}{\sqrt{\pi}}\int_{0}^{z}e^{-t^2}\, dt
\end{equation}
and
\begin{equation*}
\textup{erfc}(z)=1-\textup{erf}(z)=\frac{2}{\sqrt{\pi}}\int_{z}^{\infty}e^{-t^2}\, dt
\end{equation*}
respectively, are two important special functions having a number of applications in probability theory, statistics, physics and partial differential equations. In probability, they are related to the Gaussian normal distribution. Glashier \cite{glashier1} was the first person to coin the term \emph{Error-function} and then the term \emph{Error-function complement} in the sequel \cite{glashier2}. However, his definitions are exactly opposite to the standard definitions given above and do not involve the normalization factor $2/\sqrt{\pi}$. 
%In these two papers, Glashier evaluates a number of series and integrals involving the error functions. He also discusses the work by Kramp from his memoir \cite{kramp} which essentially involves the error functions and mentions their applications in the theory of heat conduction and in astronomical and terrestrial refraction.
 
The imaginary error function $\textup{erfi}(z)$ is defined by \cite[p.~32]{jl}\footnote{This definition differs from a factor of $\frac{2}{\sqrt{\pi}}$ from the definition in \cite{jl}.}
\begin{equation}\label{erfiz}
\textup{erfi}(z)=\frac{2}{\sqrt{\pi}}\int_{0}^{z}e^{t^2}\, dt.
\end{equation}
From (\ref{erfz}) and (\ref{erfiz}), it is straightforward that 
\begin{equation}\label{errrel}
\textup{erf}(iz)=i\textup{erfi}(z).
\end{equation}

We now give below the transformation linking the integrals of the type in \eqref{intanpar} with the error functions $\textup{erf}(z)$ and $\textup{erfi}(z)$. This transformation is of the form $G(z,\alpha)=G(iz,\beta)$ for $\alpha,\beta >0$, $\alpha\beta=1$ and $z\in\mathbb{C}$. It is also related to an integral involving the Riemann $\Xi$-function, which is defined by
\begin{equation}\label{xif}
\Xi(t):=\xi(\tf{1}{2}+it),
\end{equation}
where $\xi(s)$ is Riemann's $\xi$-function defined by \cite[p.~60]{dav}
\begin{equation}\label{rixi}
\xi(s):=\frac{1}{2}s(s-1)\pi^{-s/2}\Gamma\left(\frac{s}{2}\right)\zeta(s),
\end{equation}
$\Gamma(s)$ and $\zeta(s)$ being the Gamma function and the Riemann zeta function respectively.

\begin{theorem}\label{genmrramg}
Let $z\in\mathbb{C}, \alpha>0$ and let $\Delta\left(\alpha,z,\frac{1+it}{2}\right)$ be defined by
\begin{equation}\label{delta}
\Delta(x,z,s):=\omega(x,z,s)+\omega(x,z,1-s),
\end{equation}
with 
\begin{align}\label{omega}
\omega(x,z,s)&:=x^{\frac{1}{2}-s}e^{-\frac{z^2}{8}}{}_1F_{1}\left(1-\frac{s}{2};\frac{3}{2};\frac{z^2}{4}\right),
\end{align}
where ${}_1F_{1}(a;c;z)$ is Kummer's confluent hypergeometric function. Let $\Xi(t)$ be defined in \eqref{xif}. Then for $\alpha\beta=1$,
\begin{align}\label{mrramg}
&\sqrt{\alpha}e^{\frac{z^2}{8}}\left(\textup{erf}\left(\frac{z}{2}\right)-4\int_{0}^{\infty}\frac{e^{-\pi\alpha^2 x^2}\sin(\sqrt{\pi}\alpha xz)}{e^{2\pi x}-1}\, dx\right)\nonumber\\
&=\sqrt{\beta}e^{\frac{-z^2}{8}}\left(\textup{erfi}\left(\frac{z}{2}\right)-4\int_{0}^{\infty}\frac{e^{-\pi\beta^2 x^2}\sinh(\sqrt{\pi}\beta xz)}{e^{2\pi x}-1}\, dx\right)\nonumber\\
&=\frac{z}{8\pi^2}\int_{0}^{\infty}\Gamma\left(\frac{-1+it}{4}\right)\Gamma\left(\frac{-1-it}{4}\right)\Xi\left(\frac{t}{2}\right)\Delta\left(\alpha,z,\frac{1+it}{2}\right)\, dt,
\end{align}
where $\textup{erf}(z)$ and $\textup{erfi}(z)$ are defined in \eqref{erfz} and \eqref{erfiz} respectively.
\end{theorem}
We prove Theorem \ref{genmrramg} by evaluating the integral on the extreme right of \eqref{mrramg} to be equal to the extreme left, and by exploiting the fact that it is invariant under the simultaneous replacement of $\a$ by $\b$ and $z$ by $iz$.

This transformation generalizes a formula of Ramanujan which he wrote in his first letter to Hardy \cite[p.~XXVI]{ramcol} and which also appears in \cite[Equation (13)]{sdi}. This formula is equivalent to the first equality in the following identity, valid for $\alpha\beta=1$, and which is also due to Ramanujan \cite{riemann}:
\begin{align}\label{mrram}
&\alpha^{\frac{1}{2}}-4\pi\alpha^{\frac{3}{2}}\int_{0}^{\infty}\frac{xe^{-\pi\alpha^2 x^2}}{e^{2\pi x}-1}\, dx=\beta^{\frac{1}{2}}-4\pi\beta^{\frac{3}{2}}\int_{0}^{\infty}\frac{xe^{-\pi\beta^2 x^2}}{e^{2\pi x}-1}\, dx\nonumber\\
&=\frac{1}{4\pi\sqrt{\pi}}\int_{0}^{\infty}\Gamma\left(\frac{-1+it}{4}\right)\Gamma\left(\frac{-1-it}{4}\right)\Xi\left(\frac{t}{2}\right)\cos \left(\frac{1}{2}t\log\alpha\right)\, dt.
\end{align}
Mordell \cite[p.~331]{mordell2} rephrased the first equality in the above formula in the form \footnote{There is a misprint in Mordell's formulation of Equation \eqref{morpara}, namely, there is an extra minus sign in front of the right-hand side which should not be present.}
\begin{equation}\label{morpara}
\int_{-\infty}^{\infty}\frac{xe^{-\pi ix^2/\tau}}{e^{2\pi x}-1}\, dx=(-i\tau)^{3/2}\int_{-\infty}^{\infty}\frac{xe^{-\pi i\tau x^2}}{e^{2\pi x}-1}\, dx.
\end{equation} 
That (\ref{mrram}) is a special case of Theorem \ref{genmrramg} is not difficult to derive: for $z\neq 0$, divide both sides by $z$, let $z\to 0$ and note that 
\begin{equation*}
\lim_{z\to 0}\frac{\textup{erf}(z)}{z}=\frac{2}{\sqrt{\pi}}=\lim_{z\to 0}\frac{\textup{erfi}(z)}{z}.
\end{equation*}
 
A one-variable generalization of the integral on the extreme right-hand side in \eqref{mrram} was given in \cite[Theorem 1.5]{dixthet}, which in turn gave a generalization of the extreme left side. However, this general integral is not invariant under the simultaneous application of $\alpha\to\beta$ and $z\to iz$, and so a transformation formula generalizing the first equality in \eqref{mrram} could not be obtained there. This shortcoming is overcome in Theorem \ref{genmrramg}.

We also obtain another transformation involving error functions that is complementary to the one in Theorem \ref{genmrramg}.
\begin{theorem}\label{thmseft}
Let $z\in\mathbb{C}$. For $\a, \b>0, \a\b=1$,
\begin{align}\label{seft2}
&\sqrt{\alpha}e^{\frac{z^2}{8}}\left(\textup{erf}\left(\frac{z}{2}\right)+4\int_{-\infty}^{0}\frac{e^{-\pi\alpha^2 x^2}\sin(\sqrt{\pi}\alpha xz)}{e^{2\pi x}-1}\, dx\right)\nonumber\\
&=\sqrt{\beta}e^{\frac{-z^2}{8}}\left(\textup{erfi}\left(\frac{z}{2}\right)+4\int_{-\infty}^{0}\frac{e^{-\pi\beta^2 x^2}\sinh(\sqrt{\pi}\beta xz)}{e^{2\pi x}-1}\, dx\right).\nonumber\\
%&=\frac{z}{8\pi^2}\int_{0}^{\infty}\Gamma\left(\frac{-1+it}{4}\right)\Gamma\left(\frac{-1-it}{4}\right)\Xi\left(\frac{t}{2}\right)\Delta\left(\alpha,z,\frac{1+it}{2}\right)\, dt,
\end{align}
\end{theorem}
It is important to observe here that subtracting the corresponding sides of the first equality in \eqref{mrramg} from those of \eqref{seft2} results in \eqref{rt}, thus providing a new proof of \eqref{rt}, and hence of \eqref{fwttrans}.

%Note that we have 
%\begin{equation}\label{erfr}
%\textup{erf}(iz)=i\textup{erfi}(z).
%\end{equation}

Let $\chi$ denote a primitive Dirichlet character modulo $q$. The character analogue $\Xi(t,\chi)$ of $\Xi(t)$ is given by
\begin{equation*}
\Xi(t,\chi):=\xi\left(\frac{1}{2}+it,\chi\right),
\end{equation*}
where $\xi(s,\chi):=\left(\pi/q\right)^{-(s+a)/2}\Gamma\left(\frac{s+a}{2}\right)L(s,\chi)$, and $a=0$ if $\chi(-1)=1$ and $a=1$ if $\chi(-1)=-1$. The functional equation of $\xi(s,\chi)$ is given by $\xi(1-s,\overline{\chi})=\epsilon (\chi)\xi(s,\chi)$, where $\epsilon (\chi)=i^{a}q^{1/2}/G(\chi)$ and $G(\chi)=\sum_{m=1}^{q}\chi(m)e^{2\pi im/q}$ is the Gauss sum. See \cite[p.~69-72]{dav}. For real primitive characters, we have
\begin{equation*}
G(\chi)=
\begin{cases}
\sqrt{q}, \quad\mbox{for}\hspace{1mm}\chi\hspace{1mm}\mbox{even},\\
i\sqrt{q}, \quad\mbox{for}\hspace{1mm}\chi\hspace{1mm}\mbox{odd}.\\
\end{cases}
\end{equation*}
Hence the functional equation in this case reduces to $\xi(1-s,\chi)=\xi(s,\chi)$, which also gives $\Xi(-t,\chi)=\Xi(t,\chi)$.

We now give below the analogues of Theorem \ref{genmrramg} for real primitive characters.
\begin{theorem}\label{thmrgenrchi}
Let $z\in\mathbb{C}$ and let $\alpha$ and $\beta$ be positive numbers such that $\alpha\beta=1$. Let $\chi$ be a real primitive Dirichlet character modulo $q$.\\ 

\textup{(i)} If $\chi$ is even, 
\begin{align}\label{rgenrchi}
&\sqrt{\alpha}e^{\frac{z^2}{8}}\int_{0}^{\infty}e^{-\frac{\pi\alpha^2 x^2}{q}}\sin\left(\frac{\sqrt{\pi}\alpha xz}{\sqrt{q}}\right)\frac{\sum_{r=1}^{q-1}\chi(r)e^{-\frac{2\pi r x}{q}}}{1-e^{-2\pi x}}\, dx\nonumber\\
&=\sqrt{\beta}e^{\frac{-z^2}{8}}\int_{0}^{\infty}e^{-\frac{\pi\beta^2 x^2}{q}}\sinh\left(\frac{\sqrt{\pi}\beta xz}{\sqrt{q}}\right)\frac{\sum_{r=1}^{q-1}\chi(r)e^{-\frac{2\pi r x}{q}}}{1-e^{-2\pi x}}\, dx\nonumber\\
&=\frac{z\sqrt{q}}{16\pi^{2}}\int_{0}^{\infty}\Gamma\left(\frac{3+it}{4}\right)\Gamma\left(\frac{3-it}{4}\right)\Xi\left(\frac{t}{2},\chi\right)\Delta\left(\alpha,z,\frac{1+it}{2}\right)\, dt.
\end{align}

\textup{(ii)} If $\chi$ is odd, 
\begin{align}\label{rgenrchio}
&\sqrt{\alpha}e^{\frac{z^2}{8}}\int_{0}^{\infty}e^{-\frac{\pi\alpha^2 x^2}{q}}\cos\left(\frac{\sqrt{\pi}\alpha xz}{\sqrt{q}}\right)\frac{\sum_{r=1}^{q-1}\chi(r)e^{-\frac{2\pi r x}{q}}}{1-e^{-2\pi x}}\, dx\nonumber\\
&=\sqrt{\beta}e^{\frac{-z^2}{8}}\int_{0}^{\infty}e^{-\frac{\pi\beta^2 x^2}{q}}\cosh\left(\frac{\sqrt{\pi}\beta xz}{\sqrt{q}}\right)\frac{\sum_{r=1}^{q-1}\chi(r)e^{-\frac{2\pi r x}{q}}}{1-e^{-2\pi x}}\, dx\nonumber\\
&=\frac{1}{16\pi^{\frac{3}{2}}}\int_{0}^{\infty}\Gamma\left(\frac{1+it}{4}\right)\Gamma\left(\frac{1-it}{4}\right)\Xi\left(\frac{t}{2},\chi\right)\nabla\left(\alpha,z,\frac{1+it}{2}\right)\, dt.
\end{align}
\end{theorem}
Note that the sums inside the integrals in the above theorem are Gauss sums of purely imaginary arguments.\\ 

The first equality in \eqref{mrram} can be rewritten for $\alpha\beta=\pi^2$ as
\begin{equation}\label{mtt}
\alpha^{-\frac{1}{4}}\left(1+4\alpha\int_{0}^{\infty}\frac{xe^{-\alpha x^2}}{e^{2\pi x}-1}\, dx\right)=\beta^{-\frac{1}{4}}\left(1+4\beta\int_{0}^{\infty}\frac{xe^{-\beta x^2}}{e^{2\pi x}-1}\, dx\right).
\end{equation}
In \cite[Volume 2, p.~268]{ramnot}, Ramanujan gives an elegant approximation to the above expressions.\\

\textit{Let $\alpha>0$, $\beta>0$, $\alpha\beta=\pi^2$. Define
\begin{equation*}
I(\alpha):=\alpha^{-\frac{1}{4}}\left(1+4\alpha\int_{0}^{\infty}\frac{xe^{-\alpha x^2}}{e^{2\pi x}-1}\, dx\right).
\end{equation*}
Then 
\begin{equation}\label{near}
I(\alpha)=\left(\frac{1}{\alpha}+\frac{1}{\beta}+\frac{2}{3}\right)^{1/4}, \text{``nearly''}.
\end{equation}}\\

As mentioned by Berndt and Evans in \cite{br}, Ramanujan frequently used the words ``nearly'' or ``very nearly'' at the end of his asymptotic expansions and approximations. The above approximation is very good for values of $\alpha$ that are either very small or very large. A proof of the above fact was given in \cite{br}, where as an intermediate result, the following asymptotic expansion for $I(\a)$ as $\a\to 0$ was first obtained:
\begin{equation*}
I(\a)\sim\frac{1}{\a^{1/4}}+\frac{\a^{3/4}}{6}-\frac{\a^{7/4}}{60}+\cdots.
\end{equation*}
%\begin{proof}
%We know that for $|x|<2\pi$,
%\begin{equation}
%\frac{x}{e^x-1}=\sum_{n=0}^{\infty}\frac{B_nx^{n}}{n!}.
%\end{equation}
 
%Thus using it, we find that as $\beta\to\infty$,
%\begin{equation}\label{awl}
%I(\beta)\sim\beta^{-1/4}+\sum_{n=0}^{\infty}I_n,
%\end{equation}
%where 
%\begin{equation}
%I_n=4\beta^{3/4}\frac{{B_n}(2\pi)^{n-1}}{n!}\int_{0}^{\infty}e^{-\beta x^2}x^n\, dx.
%\end{equation}
%Elementary calculations show that
%\begin{equation}
%I_0=\alpha^{-1/4}, I_1=-\beta^{-1/4}, I_2=\tfrac{1}{6}\alpha^{3/4},\hspace{1mm}\text{and}\hspace{1mm}I_4=-\tfrac{1}{60}\alpha^{7/4}.
%\end{equation}
%Substituting the above values in \eqref{awl}, we get the following asymptotic expansion as $\beta\to\infty$, or equivalently as $\alpha\to 0$:

%\begin{equation}\label{asyma}
%I(\alpha)=I(\beta)\sim\frac{1}{\alpha^{1/4}}+\frac{\alpha^{3/4}}{6}-\frac{\alpha^{7/4}}{60}+\cdots,
%\end{equation}
%where the first equality comes from \eqref{mtt}. 

%On the other hand, by Taylor's theorem, for $\alpha$ sufficiently small, or equivalently for $\beta$ sufficiently large,
%\begin{equation}\label{asyma1}
%\left(\frac{1}{\alpha}+\frac{1}{\beta}+\frac{2}{3}\right)^{1/4}=\frac{1}{\alpha^{1/4}}+\frac{\alpha^{3/4}}{6}-\left(\frac{1}{24}-\frac{1}{4\pi^2}\right)\alpha^{7/4}+\cdots.
%\end{equation}
%Now since $\frac{1}{60}=0.01666\cdots$ and $\frac{1}{24}-\frac{1}{4\pi^2}=0.01633\cdots$, (\ref{near}) follows.
%\end{proof}
Observe that for $\alpha\beta=\pi^2$ and $z\neq 0$, the first equality in Theorem \ref{genmrramg} can be rephrased as follows:
\begin{align}\label{gen}
I(z,\a)&:=\frac{\sqrt{\pi}}{z}\alpha^{-\frac{1}{4}}e^{\frac{z^2}{8}}\textup{erf}\left(\frac{z}{2}\right)+\frac{4}{z}\alpha^{\frac{1}{4}}e^{-\frac{z^2}{8}}\int_{0}^{\infty}\frac{e^{-\alpha x^2}\sinh(\sqrt{\alpha} xz)}{e^{2\pi x}-1}\, dx\nonumber\\
&=\frac{\sqrt{\pi}}{z}\beta^{-\frac{1}{4}}e^{-\frac{z^2}{8}}\textup{erfi}\left(\frac{z}{2}\right)+\frac{4}{z}\beta^{\frac{1}{4}}e^{\frac{z^2}{8}}\int_{0}^{\infty}\frac{e^{-\beta x^2}\sin(\sqrt{\beta} xz)}{e^{2\pi x}-1}\, dx=:I(iz,\b),
\end{align}
of which \eqref{mtt} is the special case when $z\to 0$. The following general asymptotic expansion holds for the two sides in the above identity as $\a\to 0$, or equivalently as $\beta\to\infty$.
\begin{theorem}\label{1asy}
Fix $z\in\mathbb{C}$. 
%\begin{equation}
%I(z,\a):=\frac{\sqrt{\pi}}{z}\alpha^{-\frac{1}{4}}e^{\frac{z^2}{8}}\textup{erf}\left(\frac{z}{2}\right)+\frac{4}{z}\alpha^{\frac{1}{4}}e^{-\frac{z^2}{8}}\int_{0}^{\infty}\frac{e^{-\alpha x^2}\sinh(\sqrt{\alpha} xz)}{e^{2\pi x}-1}\, dx.
%\end{equation}
As $\a\to 0$,
\begin{align}\label{1asyeq}
I(z,\a)&\sim\frac{-2}{\sqrt{\pi}}\a^{-1/4}e^{z^2/8}\sum_{m=0}^{\infty}\left(\frac{-\a}{\pi^2}\right)^{m}\zeta(2m)\G\left(m+\frac{1}{2}\right){}_1F_{1}\left(m+\frac{1}{2};\frac{3}{2};\frac{-z^2}{4}\right).
\end{align}
That is,
\begin{align*}
I(z,\a)&\sim\frac{\sqrt{\pi}}{z}e^{z^2/8}\textup{erf}\left(\frac{z}{2}\right)\alpha^{-1/4}+\frac{e^{-z^2/8}}{6}\alpha^{3/4}+\frac{(z^2-6)e^{-z^2/8}}{360}\alpha^{7/4}\nonumber\\
&\quad+\frac{(60-20z^2+z^4)e^{-z^2/8}}{15120}\alpha^{11/4}+\cdots.
\end{align*}
\end{theorem}
Note that both sides of \eqref{gen} are even functions of $z$. If we successively differentiate \eqref{gen} $n$ times with respect to $z$ and then let $z\to 0$, we do not get anything interesting for odd $n$. However for $n$ even, two different behaviors are noted.

First, $n\equiv0\hspace{1mm}(\text{mod}\hspace{1mm} 4)$, i.e., $n=4k, k\in\mathbb{N}\cup\{0\}$, gives the following transformation of the form $H_{k}(\a)=H_{k}(\b)$.
\begin{theorem}\label{thmtrn0m4}
Let $\a\b=\pi^2$. Then for a non-negative integer $k$,
\begin{align}\label{trn0m4}
&\a^{-1/4}{}_2F_{1}\left(-2k,1;\frac{3}{2};2\right)+4\a^{3/4}\int_{0}^{\infty}\frac{xe^{-\a x^2}}{e^{2\pi x}-1}{}_1F_{1}\left(-2k;\frac{3}{2};2\a x^2\right)\, dx\nonumber\\
&=\b^{-1/4}{}_2F_{1}\left(-2k,1;\frac{3}{2};2\right)+4\b^{3/4}\int_{0}^{\infty}\frac{xe^{-\b x^2}}{e^{2\pi x}-1}{}_1F_{1}\left(-2k;\frac{3}{2};2\b x^2\right)\, dx.
\end{align}
\end{theorem}
Ramanujan's approximation in \eqref{near} is a special case, when $k=0$, of the following result:
\begin{theorem}\label{thmgennear}
Let $k$ be a non-negative integer. Both sides of \eqref{trn0m4} are approximated by
\begin{align}\label{gennear}
&\a^{-1/4}{}_2F_{1}\left(-2k,1;\frac{3}{2};2\right)+4\a^{3/4}\int_{0}^{\infty}\frac{xe^{-\a x^2}}{e^{2\pi x}-1}{}_1F_{1}\left(-2k;\frac{3}{2};2\a x^2\right)\, dx\nonumber\\
&={}_2F_{1}\left(-2k,1;\frac{3}{2};2\right)\left(\frac{1}{\a}+\frac{1}{\b}+\frac{2}{3\cdot{}_2F_{1}\left(-2k,1;\frac{3}{2};2\right)}\right)^{1/4}, \text{``nearly''}.
\end{align}
\end{theorem}
Again, the above right side is a very good approximation of the left side for the values of $\a$ that are either very small or very large. 
%When $\a=\b=\pi$, in particular, we obtain the following corollary.
%\begin{corollary}
%\begin{align}\label{approx}
%&\int_{0}^{\infty}\frac{e^{-y}}{e^{2\sqrt{\pi y}}-1}{}_1F_{1}\left(-2k;\frac{3}{2};2y\right)\, dy\nonumber\\
%&=\frac{1}{2}{}_2F_{1}\left(-2k,1;\frac{3}{2};2\right)\left\{\left(2+\frac{2\pi}{3\cdot{}_2F_{1}\left(-2k,1;\frac{3}{2};2\right)}\right)^{1/4}-1\right\}, \text{``nearly''}.
%\end{align}
%\end{corollary}
%Table 1 compares the two sides of the above formula for some specific values of $k$.

When $n\equiv2\hspace{1mm}(\text{mod}\hspace{1mm} 4)$, i.e., $n=4k+2, k\in\mathbb{N}\cup\{0\}$, we get a transformation of the form $J_{k}(\a)=-J_{k}(\b)$ given below. 
\begin{theorem}\label{thmtrn2m4}
Let $\a\b=\pi^2$ and let $k$ be a non-negative integer. Then,
\begin{align}\label{trn2m4}
&J_{k}(\a):=\a^{-1/4}{}_2F_{1}\left(-2k-1,1;\frac{3}{2};2\right)+4\a^{3/4}\int_{0}^{\infty}\frac{xe^{-\a x^2}}{e^{2\pi x}-1}{}_1F_{1}\left(-2k-1;\frac{3}{2};2\a x^2\right)\, dx\nonumber\\
&=-\b^{-1/4}{}_2F_{1}\left(-2k-1,1;\frac{3}{2};2\right)-4\b^{3/4}\int_{0}^{\infty}\frac{xe^{-\b x^2}}{e^{2\pi x}-1}{}_1F_{1}\left(-2k-1;\frac{3}{2};2\b x^2\right)dx=:-J_{k}(\b).
\end{align}
\end{theorem}
In particular, $J_{k}(\pi)=0$, which results in a beautiful exact evaluation of the integral in \eqref{trn2m4}.
\begin{corollary}\label{corexact}
For any non-negative integer $k$,
\begin{align}\label{exact}
\int_{0}^{\infty}\frac{xe^{-\pi x^2}}{e^{2\pi x}-1}{}_1F_{1}\left(-2k-1;\frac{3}{2};2\pi x^2\right)\, dx=-\frac{1}{4\pi}{}_2F_{1}\left(-2k-1,1;\frac{3}{2};2\right).
\end{align}
\end{corollary}
Results corresponding to the ones in Theorem \ref{thmtrn0m4} - Corollary \ref{corexact} that can be obtained by writing Theorem \ref{thmseft} in an alternate form (see \eqref{sefty} below) are collectively put in Theorem \ref{appseft} at the end of Section 5. When we combine the results from Theorem \ref{thmtrn0m4}-Corollary \ref{corexact} with those in Theorem \ref{appseft}, we obtain the following interesting theorem.
\begin{theorem}\label{ramtran}
Let $\a, \b$ be two positive numbers such that $\a\b=\pi^2$ and let $k$ be any non-negative integer. Then
{\allowdisplaybreaks\begin{align}
&\textup{(i)}\hspace{2mm} \a^{3/4}\int_{-\infty}^{\infty}\frac{xe^{-\a x^2}}{e^{2\pi x}-1}{}_1F_{1}\left(-2k;\frac{3}{2};2\a x^2\right)\, dx=\b^{3/4}\int_{-\infty}^{\infty}\frac{xe^{-\b x^2}}{e^{2\pi x}-1}{}_1F_{1}\left(-2k;\frac{3}{2};2\b x^2\right)\, dx.\nonumber\\
&\textup{(ii)}\hspace{2mm} \a^{3/4}\int_{-\infty}^{\infty}\frac{xe^{-\a x^2}}{e^{2\pi x}-1}{}_1F_{1}\left(-2k;\frac{3}{2};2\a x^2\right)\, dx\nonumber\\
&\quad\quad=\frac{1}{2}\cdot{}_2F_{1}\left(-2k,1;\frac{3}{2};2\right)\left(\frac{1}{\a}+\frac{1}{\b}+\frac{2}{3\cdot{}_2F_{1}\left(-2k,1;\frac{3}{2};2\right)}\right)^{1/4}, \text{``nearly''}.\label{gennearg}\\
&\textup{(iii)}\hspace{2mm} \a^{3/4}\int_{-\infty}^{\infty}\frac{xe^{-\a x^2}}{e^{2\pi x}-1}{}_1F_{1}\left(-2k-1;\frac{3}{2};2\a x^2\right)\, dx\nonumber\\
&\quad\quad\quad\quad\quad\quad\quad\quad\quad=-\b^{3/4}\int_{-\infty}^{\infty}\frac{xe^{-\b x^2}}{e^{2\pi x}-1}{}_1F_{1}\left(-2k-1;\frac{3}{2};2\b x^2\right)\, dx.\label{dires}
\end{align}}
In particular when $\a=\b=\pi$, we have
\begin{equation}\label{exact0}
\int_{-\infty}^{\infty}\frac{xe^{-\pi x^2}}{e^{2\pi x}-1}{}_1F_{1}\left(-2k-1;\frac{3}{2};2\pi x^2\right)\, dx=0.
\end{equation}
\end{theorem}

In \cite{riemann}, Ramanujan considered two integrals, one being that on the extreme right of \eqref{mrram}, and the second one given by
\begin{multline}
\int_{0}^{\infty} 
\Gamma \left( \frac{z-1+ i t }{4} \right)
\Gamma \left( \frac{z-1- i t }{4} \right)\Xi \left( \frac{t + i z}{2} \right) 
\Xi \left( \frac{t - i z}{2} \right) 
\frac{\cos( \tfrac{1}{2} t \log \alpha) \, dt}{t^{2} + (z+1)^{2}}.
\label{ramanujan-98}
\end{multline}
Oloa \cite[Equation 1.5]{oloa} found the asymptotic expansion\footnote{There is a misprint in this asymptotic expansion given in Oloa's paper. The minus sign in front of the second expression on the right should be a plus.} of the special case of this integral when $z=0$, namely, as $\a\to\infty$,

\textit{\begin{align*}
&\frac{1}{\pi^{3/2}}\int_0^{\infty}\Xi^{2}\left(\frac{1}{2}t\right)\left|\Gamma\left(\frac{-1+it}{4}\right)\right|^2
\frac{\cos\left(\frac{1}{2}t\log\alpha\right)}{1+t^2}\, dt\\
&\sim\frac{1}{2}\frac{\log\alpha}{\sqrt{\alpha}}+\frac{1}{2\sqrt{\alpha}}\left(\log 2\pi-\gamma\right)+\frac{\pi^{2}}{72\alpha^{3/2}}-\frac{\pi^4}{10800\alpha^{7/2}}+\cdots.
\end{align*}}

In the following theorem, we obtain the asymptotic expansion of the general integral \eqref{ramanujan-98} as $\a\to\infty$.
\begin{theorem}\label{2asy}
Fix $z$ such that $-1 <$ \textup{Re} $z < 1$. As $\alpha\to\infty$,
\begin{align}\label{oloag}
&\frac{1}{\pi^{\frac{z+3}{2}}}\int_{0}^{\infty} 
\Gamma \left( \frac{z-1+it}{4} \right)
\Gamma \left( \frac{z-1-it}{4} \right)
\Xi \left( \frac{t + iz}{2} \right)
\Xi \left( \frac{t - iz}{2} \right)\frac{\cos( \tfrac{1}{2} t \log \alpha) \, dt}{t^{2} + (z+1)^{2}}\nonumber\\
&\sim-\frac{\Gamma(z)\zeta(z)\a^{\frac{z-1}{2}}}{(2\pi)^z}-\frac{\G(z+1)\zeta(z+1)}{2\a^{\frac{z+1}{2}}(2\pi)^z}\nonumber\\
&\quad+2\a^{\frac{z+1}{2}}\sum_{m=0}^{\infty}\frac{(-1)^m}{(2\pi\a)^{2m+z+2}}\G(2m+2+z)\zeta(2m+2)\zeta(2m+z+2).
\end{align}
\end{theorem}
%Equation (1.6) from our paper gives for $\a\b=1$,
%\begin{align}\label{mrram}
%&\alpha^{\frac{1}{2}}-4\pi\alpha^{\frac{3}{2}}\int_{0}^{\infty}\frac{xe^{-\pi\alpha^2 x^2}}{e^{2\pi x}-1}\, dx=\beta^{\frac{1}{2}}-4\pi\beta^{\frac{3}{2}}\int_{0}^{\infty}\frac{xe^{-\pi\beta^2 x^2}}{e^{2\pi x}-1}\, dx\nonumber\\
%&=\frac{1}{4\pi\sqrt{\pi}}\int_{0}^{\infty}\Gamma\left(\frac{-1+it}{4}\right)\Gamma\left(\frac{-1-it}{4}\right)\Xi\left(\frac{t}{2}\right)\cos \left(\frac{1}{2}t\log\alpha\right)\, dt.
%\end{align}
This paper is organized as follows. In Section 2, we state preliminary theorems and lemmas that are subsequently used. Section 3 contains proofs of Theorems \ref{genmrramg} and \ref{thmrgenrchi}. In Section 4, we prove Theorems \ref{thmseft} and \ref{1asy}. The analogues of Theorem \ref{1asy} corresponding to the second error function transformation and to Ramanujan's transformation \eqref{rt} are also given in this section. Section 5 is devoted to proofs of Theorem \ref{thmtrn0m4} - Corollary \ref{corexact} and their analogues. We prove Theorem \ref{2asy} in Section 6. Finally, Section 7 is reserved for some concluding remarks and open problems.
\section{Nuts and bolts}
Let $f$ be an even function of $t$ of the form $f(t)=\phi(it)\phi(-it)$, where $\phi$ is analytic in $t$ as a function of a real variable. Using the functional equation for $\zeta(s)$ in the form $\xi(s)=\xi(1-s)$, where $\xi(s)$ is defined in \eqref{rixi}, it is easy to obtain the following line integral representation for the integral on the left side below, of which the integral on the extreme right of \eqref{mrramg} is a special case
\begin{equation}\label{sp4}
\int_{0}^{\infty}f\left(\frac{t}{2}\right)\Xi\left(\frac{t}{2}\right)\Delta\left(\alpha,z,\frac{1+it}{2}\right)\, dt
=\frac{2}{i}\int_{\frac{1}{2}-i\infty}^{\frac{1}{2}+i\infty}\phi\left(s-\frac{1}{2}\right)\phi\left(\frac{1}{2}-s\right)\xi(s)\omega(\alpha,z,s)\, ds,
\end{equation}
whenever the integral on the left converges. Here $\Delta(x, z, s)$ and $\omega(x, z, s)$ are the same as defined in \eqref{delta} and \eqref{omega}. Analogous to these, define
\begin{equation}\label{nabla}
\nabla(x,z,s):=\rho(x,z,s)+\rho(x,z,1-s),
\end{equation}
where 
\begin{align*}
\rho(x,z,s):=x^{\frac{1}{2}-s}e^{-\frac{z^2}{8}}{}_1F_{1}\left(\frac{1-s}{2};\frac{1}{2};\frac{z^2}{4}\right).
\end{align*}
Then for $\chi$, a real primitive character modulo $q$, the following formulas can be similarly proved.
\begin{align}\label{sp4chi}
\int_{0}^{\infty}f\left(\frac{t}{2}\right)\Xi\left(\frac{t}{2},\chi\right)\nabla\left(\alpha,z,\frac{1+it}{2}\right)\, dt
&=\frac{2}{i}\int_{\frac{1}{2}-i\infty}^{\frac{1}{2}+i\infty}\phi\left(s-\frac{1}{2}\right)\phi\left(\frac{1}{2}-s\right)\xi(s,\chi)\rho(\alpha,z,s)\, ds,\nonumber\\
\int_{0}^{\infty}f\left(\frac{t}{2}\right)\Xi\left(\frac{t}{2},\chi\right)\Delta\left(\alpha,z,\frac{1+it}{2}\right)\, dt
&=\frac{2}{i}\int_{\frac{1}{2}-i\infty}^{\frac{1}{2}+i\infty}\phi\left(s-\frac{1}{2}\right)\phi\left(\frac{1}{2}-s\right)\xi(s,\chi)\omega(\alpha,z,s)\, ds,
\end{align}
whenever the integrals on the left-hand sides converge. Note that for $\alpha\beta=1$,
\begin{equation}\label{rel}
\Delta\left(\alpha,z,\frac{1+it}{2}\right)=\Delta\left(\beta,iz,\frac{1+it}{2}\right), \nabla\left(\alpha,z,\frac{1+it}{2}\right)=\nabla\left(\beta,iz,\frac{1+it}{2}\right),
\end{equation}
both of which can be proved using Kummer's first transformation for ${}_1F_{1}(a;c;w)$ \cite[p.~125, Equation (2)]{rain} given by 
\begin{equation}\label{kft}
{}_1F_{1}(a;c;w)=e^w{}_1F_{1}(c-a;c;-w).
\end{equation}
The formulas in \eqref{rel} render the integrals on the left-hand sides of \eqref{sp4} and \eqref{sp4chi} invariant under the simultaneous replacement of $\alpha$ by $\beta$ and $z$ by $iz$, and hence, as a by-product of the evaluation of these integrals, we obtain identities of the form $G(z,\alpha)=G(iz,\beta)$ and $G(z,\alpha, \chi)=G(iz,\beta,\chi)$.

In proving Theorem \ref{thmrgenrchi}, we make use of the following special case \cite[Theorem 2.1]{bds} of a result due to Berndt \cite[Theorem 10.1]{bcbspfunc}:
\begin{theorem}\label{wat}
Let $x>0$. If $\chi$ is even with period $q$ and \textup{Re} $\nu\geq 0$, then
\begin{equation*}
\sum_{n=1}^{\infty}\chi(n) n^{\nu} K_{\nu}\left(\frac{2\pi nx}{q}\right)
=\frac{\pi^{\frac{1}{2}}}{2xG(\overline{\chi})}\left(\frac{qx}{\pi}\right)^{\nu+1}\Gamma\left(\nu+\frac{1}{2}\right)
\sum_{n=1}^{\infty}\overline{\chi}(n)(n^2+x^2)^{-\nu-\frac{1}{2}};
\end{equation*}
if $\chi$ is odd with period $k$ and \textup{Re} $\nu>-1$, then
\begin{equation}\label{watodd}
\sum_{n=1}^{\infty}\chi(n) n^{\nu+1} K_{\nu}\left(\frac{2\pi
nx}{q}\right)=\frac{i\pi^{\frac{1}{2}}}{2x^2{G(\overline{\chi})}}\left(\frac{qx}{\pi}\right)^{\nu+2}
\Gamma\left(\nu+\frac{3}{2}\right)\sum_{n=1}^{\infty}\overline{\chi}(n)n(n^2+x^2)^{-\nu-\frac{3}{2}}.
\end{equation}
\end{theorem}
The following two lemmas, given in \cite[p.~503, Formula (3.952.7)]{grn} and \cite[pp.~318, 320, Formulas (10), (30)]{tti} respectively, will also be employed in the proof of Theorem \ref{thmrgenrchi}. 
\begin{lemma}\label{corinvsinl}
For $c=$ Re $s>-1$ and Re $a>0$, we have
\begin{equation*}
\frac{1}{2\pi i}\int_{c-i\infty}^{c+i\infty}\frac{b}{2}a^{-\frac{1}{2}-\frac{s}{2}}e^{-\frac{b^2}{4a}}\Gamma\left(\frac{s+1}{2}\right){}_1F_{1}\left(1-\frac{s}{2};\frac{3}{2};\frac{b^2}{4a}\right)x^{-s}\, ds=e^{-ax^2}\sin bx.
\end{equation*}
\end{lemma}
\begin{lemma}\label{corinvcosl}
For $c=$ Re $s>0$ and Re $a>0$, we have
\begin{equation}\label{invmel}
\frac{1}{2\pi i}\int_{c-i\infty}^{c+i\infty}\frac{1}{2}a^{-\tfrac{s}{2}}\Gamma\left(\frac{s}{2}\right)e^{-\frac{b^2}{4a}}{}_1F_{1}\left(\frac{1-s}{2};\frac{1}{2};\frac{b^2}{4a}\right)x^{-s}\, ds=e^{-ax^2}\cos bx.
\end{equation}
\end{lemma}
\noindent
We note that \cite[Equation (2.10)]{dixthet}
\begin{equation*}
{}_1F_{1}\left(\tfrac{1}{4}-\lambda;\tfrac{1}{2};\tfrac{z^{2}}{4}\right)\sim e^{z^2/8}\cos\left(\sqrt{\lambda}z\right),
\end{equation*}
as $|\lambda|\to\infty$ and $|\arg(\lambda z)|<2\pi$, and the Stirling's formula for $\Gamma(s)$, $s=\sigma+it$, in a vertical strip $\alpha\leq\sigma\leq\beta$ given by
\begin{equation}\label{strivert}
|\Gamma(s)|=(2\pi)^{\tf{1}{2}}|t|^{\sigma-\tf{1}{2}}e^{-\tf{1}{2}\pi |t|}\left(1+O\left(\frac{1}{|t|}\right)\right),
\end{equation}
as $|t|\to\infty$ give convergence of the integrals on the extreme right-hand sides of \eqref{mrramg}, \eqref{rgenrchi} and \eqref{rgenrchio}. If $F(s)$ and $G(s)$ denote the Mellin transforms of $f(x)$ and $g(x)$ respectively and $s$ with Re $s=c$ lies in a common strip where both $F$ and $G$ are analytic, then a variant of Parseval's formula \cite[p.83, Equation (3.1.13)]{kp} gives
\begin{equation}\label{melconv}
\frac{1}{2\pi i}\int_{c-i\infty}^{c+i\infty}F(s)G(s)w^{-s}\, ds=\int_{0}^{\infty}f(x)g\left(\frac{w}{x}\right)\frac{dx}{x}.
\end{equation}
Watson's lemma \cite[p.~71]{olver} is given by
\begin{lemma}\label{watlem}
If $q(t)$ is a function of the positive real variable $t$ such that
\begin{equation*}
q(t)\sim\sum_{s=0}^{\infty}a_st^{(s+\lambda-\mu)/\mu}\hspace{3mm} (t\to 0)
\end{equation*}
for positive constants $\lambda$ and $\mu$, then
\begin{equation}\label{watlem1}
\int_{0}^{\infty}e^{-xt}q(t)\, dt\sim\sum_{s=0}^{\infty}\Gamma\left(\frac{s+\lambda}{\mu}\right)\frac{a_s}{x^{(s+\lambda)/\mu}}\hspace{3mm} (x\to\infty),
\end{equation}
provided that this integral converges throughout its range for all sufficiently large $x$.
\end{lemma}
The above result also holds \cite[p.~32]{temme} for complex $\lambda$ with Re $\lambda>0$, and for $x\in\mathbb{C}$ with the integral being convergent for all sufficiently large values of Re $x$. 
\section{The first error function transformation and its character analogues}
We begin by proving the first error function transformation given in Theorem \ref{genmrramg} and then proceed to a proof of its character analogues given in Theorem \ref{thmrgenrchi}.
\begin{proof}[Theorem \textup{\ref{genmrramg}}][]
Let $\phi(s)=\Gamma\left(\frac{-1}{4}+\frac{s}{2}\right)$ and let
\begin{equation*}
J(z,\alpha)=\int_{0}^{\infty}\Gamma\left(\frac{-1+it}{4}\right)\Gamma\left(\frac{-1-it}{4}\right)\Xi\left(\frac{t}{2}\right)\Delta\left(\alpha,z,\frac{1+it}{2}\right)\, dt.
\end{equation*}
Use \eqref{xif}, \eqref{rixi}, \eqref{sp4}, the functional equation for $\Gamma(s)$ and the reflection formula to see that
\begin{align*}
J(z,\alpha)&=\frac{2\sqrt{\alpha}e^{-\frac{z^2}{8}}}{i}\int_{\frac{1}{2}-i\infty}^{\frac{1}{2}+i\infty}\Gamma\left(\frac{s+1}{2}\right)\Gamma\left(-\frac{s}{2}\right)\Gamma\left(1+\frac{s}{2}\right)\zeta(s){}_1F_{1}\left(1-\frac{s}{2};\frac{3}{2};\frac{z^2}{4}\right)\left(\sqrt{\pi}\alpha\right)^{-s}\, ds\nonumber\\
&=-\frac{4\sqrt{\alpha}e^{-\frac{z^2}{8}}}{i}\int_{\frac{1}{2}-i\infty}^{\frac{1}{2}+i\infty}\frac{\pi}{\sin\left(\frac{1}{2}\pi s\right)}\Gamma\left(\frac{s+1}{2}\right)\zeta(s){}_1F_{1}\left(1-\frac{s}{2};\frac{3}{2};\frac{z^2}{4}\right)\left(\sqrt{\pi}\alpha\right)^{-s}\, ds.
\end{align*}
Now in order to use the series representation for $\zeta(s)$, we shift the line of integration from Re $s=\frac{1}{2}$ to Re $s=1+\delta$, where $0<\delta<1$. Consider a positively oriented rectangular contour with sides $[\frac{1}{2}+iT, \frac{1}{2}-iT], [\frac{1}{2}-iT, 1+\delta-iT], [1+\delta-iT,1+\delta+iT]$ and $[1+\delta+iT,\frac{1}{2}+iT]$, where $T$ is any positive real number. We have to consider the contribution of the pole of order $1$ of the integrand (due to $\zeta(s)$). Using the residue theorem, noting that by (\ref{strivert}) the integrals along the horizontal line segments tend to zero as $T\to\infty$, and then interchanging the order of summation and integration we have
\begin{align*}
J(z,\alpha)&=-\frac{4\sqrt{\alpha}e^{-\frac{z^2}{8}}}{i}\bigg(\sum_{n=1}^{\infty}\int_{1+\delta-i\infty}^{1+\delta+i\infty}\frac{\pi}{\sin\left(\frac{1}{2}\pi s\right)}\Gamma\left(\frac{s+1}{2}\right)\nonumber\\
&\quad\quad\quad\quad\quad\quad\quad\quad\quad\quad\quad\quad\quad\times{}_1F_{1}\left(1-\frac{s}{2};\frac{3}{2};\frac{z^2}{4}\right)\left(\sqrt{\pi}\alpha n\right)^{-s}\, ds-2\pi iL\bigg),
\end{align*}
where
\begin{equation*}
L=\lim_{s\to 1}(s-1)\zeta(s)\frac{\pi}{\sin\left(\frac{1}{2}\pi s\right)}\Gamma\left(\frac{s+1}{2}\right){}_1F_{1}\left(1-\frac{s}{2};\frac{3}{2};\frac{z^2}{4}\right)\left(\sqrt{\pi}\alpha\right)^{-s}.
\end{equation*}
It is easy to see that 
\begin{equation*}
L=\frac{\sqrt{\pi}}{\alpha}{}_1F_{1}\left(\frac{1}{2};\frac{3}{2};\frac{z^2}{4}\right).
\end{equation*}
Now using the fact \cite[p.~98]{kp} that for $0<c=$ Re $s<2$,
\begin{equation*}
\frac{1}{2\pi i}\int_{c-i\infty}^{c+i\infty}\frac{\pi}{\sin\left(\frac{1}{2}\pi s\right)}x^{-s}\, ds=\frac{2}{(1+x^2)},
\end{equation*}
combined with the special case when $b=z\neq 0$, $a=1$ of Lemma \ref{corinvsinl}, and (\ref{melconv}), we see that
\begin{align*}
J(z,\alpha)=-8\pi\sqrt{\alpha}e^{-\frac{z^2}{8}}\bigg(\frac{4e^{\frac{z^2}{4}}}{z}\sum_{n=1}^{\infty}\int_{0}^{\infty}\frac{e^{-x^2}\sin xz}{1+\left(\frac{\sqrt{\pi}\alpha n}{x}\right)^{2}}\, \frac{dx}{x}-\frac{\sqrt{\pi}}{\alpha}{}_1F_{1}\left(\frac{1}{2};\frac{3}{2};\frac{z^2}{4}\right)\bigg).
\end{align*}
Employ the change of variable $x\to\sqrt{\pi}\alpha x$ and \eqref{kft} to see that
\begin{align}\label{jzabc}
J(z,\alpha)=-8\pi\sqrt{\alpha}e^{\frac{z^2}{8}}\bigg(\frac{4e^{\frac{z^2}{4}}}{z}\sum_{n=1}^{\infty}\int_{0}^{\infty}\frac{xe^{-\pi\alpha^2x^2}\sin\left(\sqrt{\pi}\alpha xz\right)}{n^2+x^2}\, dx-\frac{\sqrt{\pi}}{\alpha}{}_1F_{1}\left(1;\frac{3}{2};\frac{-z^2}{4}\right)\bigg).
\end{align}
Now for $t\neq 0$ \cite[p.~191]{con},
\begin{equation}\label{cotid1}
\sum_{n=1}^{\infty}\frac{1}{t^2+n^2}=\frac{\pi}{t}\left(\frac{1}{e^{2\pi t}-1}-\frac{1}{2\pi t}+\frac{1}{2}\right).
\end{equation}
Thus, interchanging the order of summation and integration in (\ref{jzabc}) and then substituting (\ref{cotid1}) and simplifying, we observe that
\begin{align}\label{jzabc1}
J(z,\alpha)&=-8\pi\sqrt{\alpha}e^{\frac{z^2}{8}}\bigg(\frac{4\pi}{z}\int_{0}^{\infty}\frac{e^{-\pi\alpha^2x^2}\sin\left(\sqrt{\pi}\alpha xz\right)}{e^{2\pi x}-1}\, dx-\frac{2}{z}\int_{0}^{\infty}\frac{e^{-\pi\alpha^2x^2}\sin\left(\sqrt{\pi}\alpha xz\right)}{x}\, dx\nonumber\\
&\quad\quad\quad\quad\quad\quad\quad+\frac{2\pi}{z}\int_{0}^{\infty}e^{-\pi\alpha^2x^2}\sin\left(\sqrt{\pi}\alpha xz\right)\, dx-\frac{\sqrt{\pi}}{\alpha}{}_1F_{1}\left(1;\frac{3}{2};\frac{-z^2}{4}\right)\bigg).
\end{align}
However from \cite[p.~488, formula 3.896, no. 3]{grn}, 
\begin{equation}\label{dawson}
\int_{0}^{\infty}e^{-\pi\alpha^2x^2}\sin\left(\sqrt{\pi}\alpha xz\right)\, dx=\frac{z}{2\sqrt{\pi}\alpha}{}_1F_{1}\left(1;\frac{3}{2};\frac{-z^2}{4}\right)
\end{equation}
and by \cite[p.~503, formula 3.952, no. 7]{grn} and \cite[p.~889, formula 8.253, no. 1]{grn} (see also \cite[p.~421]{glashier2}),
\begin{equation}\label{errapp}
\int_{0}^{\infty}\frac{e^{-\pi\alpha^2x^2}\sin\left(\sqrt{\pi}\alpha xz\right)}{x}\, dx=\frac{\pi}{2}\textup{erf}\left(\frac{z}{2}\right).
\end{equation}
Thus, substituting (\ref{dawson}) and (\ref{errapp}) in (\ref{jzabc1}) and simplifying, we finally arrive at
\begin{align}\label{mrramg1}
\frac{1}{8\pi^2}J(z,\alpha)=\frac{\sqrt{\alpha}e^{\frac{z^2}{8}}}{z}\left(\textup{erf}\left(\frac{z}{2}\right)-4\int_{0}^{\infty}\frac{e^{-\pi\alpha^2 x^2}\sin(\sqrt{\pi}\alpha xz)}{e^{2\pi x}-1}\, dx\right).
\end{align}
Using \eqref{errrel}, it is clear that simultaneously replacing $\alpha$ by $\beta$ and $z$ by $iz$ in (\ref{mrramg1}) and employing \eqref{rel} and (\ref{errrel}) give (\ref{mrramg}) since $J(z,\alpha)$ is invariant.
\end{proof}
\begin{proof}[Theorem \textup{\ref{thmrgenrchi}}][]
We prove the theorem only for odd real $\chi$. The case when $\chi$ is even and real can be similarly proved. Let $\phi(s)=\Gamma\left(\frac{1}{4}+\frac{s}{2}\right)$ and let 
\begin{equation*}
P(z,\alpha,\chi)=\int_{0}^{\infty}\Gamma\left(\frac{1+it}{4}\right)\Gamma\left(\frac{1-it}{4}\right)\Xi\left(\frac{t}{2},\chi\right)\nabla\left(\alpha,z,\frac{1+it}{2}\right)\, dt.
\end{equation*}
Using the first equality in (\ref{sp4chi}), we see that
{\allowdisplaybreaks\begin{align*}
P(z,\alpha,\chi)&=\frac{2\sqrt{\alpha q}e^{-\frac{z^2}{8}}}{i\sqrt{\pi}}\int_{\frac{1}{2}-i\infty}^{\frac{1}{2}+i\infty}\Gamma\left(\frac{s}{2}\right)\Gamma\left(\frac{1-s}{2}\right)\Gamma\left(\frac{s+1}{2}\right)L(s,\chi)\nonumber\\
&\quad\quad\quad\quad\quad\quad\quad\quad\quad\quad\times {}_1F_{1}\left(\frac{1-s}{2};\frac{1}{2};\frac{z^2}{4}\right)\left(\frac{\sqrt{\pi}\alpha}{\sqrt{q}}\right)^{-s}\, ds\nonumber\\
&=\frac{2\sqrt{\alpha q}e^{-\frac{z^2}{8}}}{i\sqrt{\pi}}\int_{\frac{1}{2}-i\infty}^{\frac{1}{2}+i\infty}\frac{\pi}{\cos\left(\frac{1}{2}\pi s\right)}\Gamma\left(\frac{s}{2}\right)L(s,\chi){}_1F_{1}\left(\frac{1-s}{2};\frac{1}{2};\frac{z^2}{4}\right)\left(\frac{\sqrt{\pi}\alpha}{\sqrt{q}}\right)^{-s}\, ds,
\end{align*}}%
where in the last step, we used a different version of the reflection formula, namely, $\Gamma\left(\frac{1}{2}+z\right)\Gamma\left(\frac{1}{2}-z\right)=\frac{\pi}{\cos \pi z}$ for $z-\frac{1}{2}\notin\mathbb{Z}$. As before, shift the line of integration from Re $s=\frac{1}{2}$ to Re $s=1+\delta$, $0<\delta<2$, employ the residue theorem and take into account the contribution from the pole of order $1$ at $s=1$ of the integrand (due to $\cos\frac{1}{2}\pi s$). This gives,
\begin{align}\label{hzacbc}
P(z,\alpha,\chi)&=\frac{2\sqrt{\alpha q}e^{-\frac{z^2}{8}}}{i\sqrt{\pi}}\bigg(\sum_{n=1}^{\infty}\chi(n)\int_{1+\delta-i\infty}^{1+\delta+i\infty}\frac{\pi}{\cos\left(\frac{1}{2}\pi s\right)}\Gamma\left(\frac{s}{2}\right)\nonumber\\
&\quad\quad\quad\quad\quad\quad\quad\quad\quad\quad\quad\quad\quad\times{}_1F_{1}\left(\frac{1-s}{2};\frac{1}{2};\frac{z^2}{4}\right)\left(\frac{\sqrt{\pi}\alpha n}{\sqrt{q}}\right)^{-s}\, ds-2\pi iL_1\bigg),
\end{align}
where
\begin{equation}\label{l1}
L_1=\lim_{s\to 1}\frac{(s-1)\pi}{\cos\left(\frac{1}{2}\pi s\right)}\Gamma\left(\frac{s}{2}\right)L(s,\chi){}_1F_{1}\left(\frac{1-s}{2};\frac{1}{2};\frac{z^2}{4}\right)\left(\frac{\sqrt{\pi}\alpha}{\sqrt{q}}\right)^{-s}=-\frac{2\sqrt{q}}{\alpha}L(1,\chi).
\end{equation}
Also replacing $s$ by $(s+1)/2$, $x$ by $x^2$ in the formula \cite[p.~91, Equation (3.3.10)]{kp}
\begin{equation*}
\frac{1}{2\pi i}\int_{c-i\infty}^{c+i\infty}\frac{x^{-s}}{\sin\pi s}\, ds=\frac{1}{\pi (1+x)},
\end{equation*}
and simplifying, we see that for $-1<$ Re $s<1$,
\begin{equation*}
\frac{1}{2\pi i}\int_{c-i\infty}^{c+i\infty}\frac{\pi}{\cos\left(\frac{1}{2}\pi s\right)}x^{-s}\, ds=\frac{2x}{1+x^2}.
\end{equation*}
Another application of the residue theorem yields for $0<c<1$, 
\begin{align}\label{anotheres}
&\int_{1+\delta-i\infty}^{1+\delta+i\infty}\frac{\pi}{\cos\left(\frac{1}{2}\pi s\right)}\Gamma\left(\frac{s}{2}\right){}_1F_{1}\left(\frac{1-s}{2};\frac{1}{2};\frac{z^2}{4}\right)\left(\frac{\sqrt{\pi}\alpha n}{\sqrt{q}}\right)^{-s}\, ds\nonumber\\
&=\int_{c-i\infty}^{c+i\infty}\frac{\pi}{\cos\left(\frac{1}{2}\pi s\right)}\Gamma\left(\frac{s}{2}\right){}_1F_{1}\left(\frac{1-s}{2};\frac{1}{2};\frac{z^2}{4}\right)\left(\frac{\sqrt{\pi}\alpha n}{\sqrt{q}}\right)^{-s}\, ds-\frac{4\pi i\sqrt{q}}{\alpha n}\nonumber\\
&=2\pi i\left(4e^{\frac{z^2}{4}}\int_{0}^{\infty}\frac{n}{x^2+n^2}e^{-\frac{\pi\alpha^2x^2}{q}}\cos\left(\frac{\sqrt{\pi}\alpha xz}{\sqrt{q}}\right)\, dx-\frac{2\sqrt{q}}{\alpha n}\right),
\end{align}
where in the last step we used Lemma \ref{corinvcosl} with $a=1$, $x=\frac{\sqrt{\pi}\alpha n}{\sqrt{q}}$ and $b=z$, and (\ref{melconv}), followed by a change of variable $x\to\frac{\sqrt{\pi}\alpha x}{\sqrt{q}}$. Now substitute (\ref{anotheres}) and (\ref{l1}) in (\ref{hzacbc}) and simplify to obtain
\begin{equation}\label{hzacbc1}
P(z,\alpha,\chi)=16\sqrt{\pi\alpha q}e^{\frac{z^2}{8}}\int_{0}^{\infty}\left(\sum_{n=1}^{\infty}\frac{n\chi(n)}{x^2+n^2}\right)e^{-\frac{\pi\alpha^2x^2}{q}}\cos\left(\frac{\sqrt{\pi}\alpha xz}{\sqrt{q}}\right)\, dx.
\end{equation}
Now use (\ref{watodd}) with $\chi$ real and $\nu=-1/2$ to see that
\begin{equation}\label{watodd2}
\sum_{n=1}^{\infty}\frac{n\chi(n)}{x^2+n^2}=\frac{\pi}{\sqrt{q}}\sum_{n=1}^{\infty}\chi(n)e^{-\frac{2\pi nx}{q}}.
\end{equation}
Employing (\ref{watodd2}) in (\ref{hzacbc1}), we have
\begin{equation}\label{hzacbc2}
P(z,\alpha,\chi)=16\sqrt{\pi^{3}\alpha}e^{\frac{z^2}{8}}\int_{0}^{\infty}e^{-\frac{\pi\alpha^2x^2}{q}}\cos\left(\frac{\sqrt{\pi}\alpha xz}{\sqrt{q}}\right)\sum_{n=1}^{\infty}\chi(n)e^{-\frac{2\pi nx}{q}}\, dx.
\end{equation}
Writing $n=mq+r, 0\leq m<\infty, 0\leq r\leq q-1$, and noting that $\chi$ is periodic with period $q$, we have
\begin{equation}\label{rc}
\sum_{n=1}^{\infty}\chi(n)e^{-\frac{2\pi nx}{q}}=\frac{\sum_{r=0}^{q-1}\chi(r)e^{-\frac{2\pi r x}{q}}}{1-e^{-2\pi x}}.
\end{equation}
Finally, (\ref{rc}) along with (\ref{hzacbc2}) gives
\begin{equation*}
\frac{1}{16\pi^{\frac{3}{2}}}P(z,\alpha,\chi)=\sqrt{\alpha}e^{\frac{z^2}{8}}\int_{0}^{\infty}e^{-\frac{\pi\alpha^2x^2}{q}}\cos\left(\frac{\sqrt{\pi}\alpha xz}{\sqrt{q}}\right)\frac{\sum_{r=1}^{q-1}\chi(r)e^{-\frac{2\pi r x}{q}}}{1-e^{-2\pi x}}\, dx.
\end{equation*}
This gives (\ref{rgenrchio}) as $P(z,\alpha,\chi)$ is invariant under the simultaneous application of the maps $\alpha\to\beta$ and $z\to iz$, which can be seen from \eqref{rel}.
\end{proof}

\section{The second error function transformation and an asymptotic expansion}
We first establish the second error function transformation given in Theorem \ref{thmseft} and then the asymptotic expansion from Theorem \ref{1asy}.
\begin{proof}[Theorem \textup{\ref{thmseft}}][]
Note that from \eqref{dawson} and \cite[p.~889, formula 8.253, no. 1]{grn},
\begin{equation}\label{fre}
\int_{0}^{\infty}e^{-\pi\alpha^2x^2}\sin\left(\sqrt{\pi}\alpha xz\right)\, dx=\frac{1}{2\a}e^{-\frac{z^2}{4}}\textup{erfi}\left(\frac{z}{2}\right).
\end{equation}
Also, 
\begin{align*}
&\sqrt{\alpha}e^{\frac{z^2}{8}}\left(\textup{erf}\left(\frac{z}{2}\right)+4\int_{-\infty}^{0}\frac{e^{-\pi\alpha^2 x^2}\sin(\sqrt{\pi}\alpha xz)}{e^{2\pi x}-1}\, dx\right)\nonumber\\
&=\sqrt{\alpha}e^{\frac{z^2}{8}}\left(\textup{erf}\left(\frac{z}{2}\right)+4\int_{0}^{\infty}e^{-\pi\alpha^2x^2}\sin\left(\sqrt{\pi}\alpha xz\right)\, dx+4\int_{0}^{\infty}\frac{e^{-\pi\alpha^2 x^2}\sin(\sqrt{\pi}\alpha xz)}{e^{2\pi x}-1}\, dx\right),
\end{align*}
where in the first step, we replaced $x$ by $-x$ in the integral, and then simplified it using $e^{2\pi x}/(e^{2\pi x}-1)=1+1/(e^{2\pi x}-1)$. Now use \eqref{fre} and the first error function transformation in \eqref{mrramg} to replace the second integral in the above equation to obtain
\begin{align}\label{seft22}
&\sqrt{\alpha}e^{\frac{z^2}{8}}\left(\textup{erf}\left(\frac{z}{2}\right)+4\int_{-\infty}^{0}\frac{e^{-\pi\alpha^2 x^2}\sin(\sqrt{\pi}\alpha xz)}{e^{2\pi x}-1}\, dx\right)\nonumber\\
&=\sqrt{\alpha}e^{\frac{z^2}{8}}\left(2\textup{erf}\left(\frac{z}{2}\right)+\frac{1}{\a}e^{-\frac{z^2}{4}}\textup{erfi}\left(\frac{z}{2}\right)+\frac{4}{\a}e^{-\frac{z^2}{4}}\int_{0}^{\infty}\frac{e^{-\pi\b^2 x^2}\sinh(\sqrt{\pi}\b xz)}{e^{2\pi x}-1}\, dx\right)\nonumber\\
&=\sqrt{\b}e^{\frac{-z^2}{8}}\left(\textup{erfi}\left(\frac{z}{2}\right)+\frac{2}{\b}e^{\frac{z^2}{4}}\textup{erf}\left(\frac{z}{2}\right)+4\int_{0}^{\infty}\frac{e^{-\pi\b^2 x^2}\sinh(\sqrt{\pi}\b xz)}{e^{2\pi x}-1}\, dx\right),
\end{align}
where we used the fact $\a\b=1$ to simplify the last step. Now replace $\a$ by $\b$ and $z$ by $iz$ in \eqref{fre}, and use \eqref{errrel} to obtain
\begin{equation}\label{fre2}
\int_{0}^{\infty}e^{-\pi\b^2x^2}\sinh\left(\sqrt{\pi}\b xz\right)\, dx=\frac{1}{2\b}e^{\frac{z^2}{4}}\textup{erf}\left(\frac{z}{2}\right).
\end{equation}
Finally, use \eqref{fre2} to simplify the extreme right of \eqref{seft22}, thereby obtaining \eqref{seft2} and thus completing the proof.
\end{proof}
\begin{proof}[Theorem \textup{\ref{1asy}}][]
By a change of variable $x^2=t$,
\begin{equation*}
\int_{0}^{\infty}\frac{e^{-\beta x^2}\sin(\sqrt{\beta} xz)}{e^{2\pi x}-1}\, dx=\int_{0}^{\infty}\frac{e^{-\beta t}\sin(z\sqrt{\beta t})}{2\sqrt{t}(e^{2\pi\sqrt{t}}-1)}\, dt.
\end{equation*}
Let 
\begin{equation*}
f(t, z)=\frac{\sin(z\sqrt{\beta t})}{2\sqrt{t}(e^{2\pi\sqrt{t}}-1)}.
\end{equation*}
First consider, for $|t|<1$,
{\allowdisplaybreaks\begin{align*}
\frac{2\pi t\sin(at)}{e^{2\pi t}-1}&=\sum_{m=0}^{\infty}\frac{B_m(2\pi)^mt^m}{m!}\sum_{n=0}^{\infty}\frac{(-1)^na^{2n+1}t^{2n+1}}{(2n+1)!}\nonumber\\
&=\sum_{j=1}^{\infty}\left(\sum_{m+2n+1=j}\frac{B_m(2\pi)^m}{m!}\frac{(-1)^na^{2n+1}}{(2n+1)!}\right)t^j\nonumber\\
&=\sum_{j=1}^{\infty}\left(\sum_{k=0\atop{k\hspace{0.5mm}\text{even}}}^{j-1}\frac{B_{j-1-k}(2\pi)^{j-1-k}}{(j-1-k)!}\frac{(-1)^{k/2}a^{k+1}}{(k+1)!}\right)t^j\nonumber\\
&=t\sum_{j=0}^{\infty}\left(\sum_{k=0\atop{k\hspace{0.5mm}\text{even}}}^{j}\frac{B_{j-k}(2\pi)^{j-k}}{(j-k)!}\frac{(-1)^{k/2}a^{k+1}}{(k+1)!}\right)t^j,
\end{align*}}
%so that
%\begin{equation}
%\frac{\sin(at)}{e^{2\pi t}-1}=\frac{1}{2\pi}\sum_{j=0}^{\infty}\left(\sum_{k=0\atop{k\hspace{0.5mm}\text{even}}}^{j}\frac{B_{j-k}(2\pi)^{j-k}}{(j-k)!}\frac{(-1)^{k/2}a^{k+1}}{(k+1)!}\right)t^j.
%\end{equation}
Replacing $t$ by $\sqrt{t}$ and $a=z\sqrt{\beta}$, we have for $|t|<1$,
\begin{equation*}
\frac{\sin(z\sqrt{\beta t})}{2\sqrt{t}(e^{2\pi\sqrt{t}}-1)}=\frac{1}{4\pi}\sum_{j=0}^{\infty}\left(\sum_{k=0\atop{k\hspace{0.5mm}\text{even}}}^{j}\frac{B_{j-k}(2\pi)^{j-k}}{(j-k)!}\frac{(-1)^{k/2}(z\sqrt{\beta})^{k+1}}{(k+1)!}\right)t^{\frac{j-1}{2}}
\end{equation*}
Thus, as $t\to 0+$,
\begin{equation*}
f(t, z)\sim\frac{1}{4\pi}\sum_{j=0}^{\infty}\left(\sum_{k=0\atop{k\hspace{0.5mm}\text{even}}}^{j}\frac{B_{j-k}(2\pi)^{j-k}}{(j-k)!}\frac{(-1)^{k/2}(z\sqrt{\beta})^{k+1}}{(k+1)!}\right)t^{\frac{j-1}{2}}
\end{equation*}
Hence by Watson's lemma, as $\beta\to\infty$,
\begin{equation}\label{watappl}
\int_{0}^{\infty}\frac{e^{-\beta t}\sin(z\sqrt{\beta t})}{2\sqrt{t}(e^{2\pi\sqrt{t}}-1)}\, dt\sim \sum_{j=0}^{\infty}\frac{\Gamma\left(\frac{j+1}{2}\right)}{\beta^{\frac{j+1}{2}}}\sum_{k=0\atop{k\hspace{0.5mm}\text{even}}}^{j}\frac{B_{j-k}(2\pi)^{j-k-1}}{2(j-k)!}\frac{(-1)^{k/2}(z\sqrt{\beta})^{k+1}}{(k+1)!}.
\end{equation}
From \eqref{watappl} and the notation in \eqref{gen}, we find that
\begin{equation}\label{ibetaz}
I(iz, \beta)\sim\frac{\sqrt{\pi}}{z}\beta^{-\frac{1}{4}}e^{-\frac{z^2}{8}}\textup{erfi}\left(\frac{z}{2}\right)+\sum_{j=0}^{\infty}I_{j,z},
\end{equation}
where
\begin{align*}
I_{j, z}=\frac{4}{z}\beta^{\frac{1}{4}}e^{\frac{z^2}{8}}\frac{\Gamma\left(\frac{j+1}{2}\right)}{\beta^{\frac{j+1}{2}}}\sum_{k=0\atop{k\hspace{0.5mm}\text{even}}}^{j}\frac{B_{j-k}(2\pi)^{j-k-1}}{2(j-k)!}\frac{(-1)^{k/2}(z\sqrt{\beta})^{k+1}}{(k+1)!}.
\end{align*}
We first evaluate $I_{j, z}$ when $j$ is odd, say $j=2n+1$. Since all of the odd-indexed Bernoulli numbers except $B_1$ are equal to zero and $B_1=-1/2$, only the last term, namely $j=2n+1$ and $k=2n$, contributes to the sum giving
\begin{equation}\label{a14}
I_{2n+1,z}=\frac{(-1)^{n+1}n!}{(2n+1)!}z^{2n}e^{z^2/8}\beta^{-1/4}.
\end{equation}
%The corresponding expression when $j$ is even is complicated since all of the terms contribute. We give below first few values of $I_{2n,z}$. For example,
%\begin{align}
%I_{0,z}&=\frac{e^{z^2/8}\beta^{1/4}}{\sqrt{\pi}}\nonumber\\
%I_{1, z}&=\frac{-ze^{z^2/8}\beta^{-1/4}}{\sqrt{\pi}}\nonumber\\
%I_{2, z}&=\frac{\sqrt{\pi}e^{z^2/8}\beta^{-3/4}}{12}\left(2\pi-\frac{z^2\beta}{\pi}\right)\nonumber\\
%I_{3, z}&=\frac{z^3e^{z^2/8}\beta^{-1/4}}{6\sqrt{\pi}}\nonumber\\
%I_{4,z}&=\frac{\sqrt{\pi}e^{z^2/8}\beta^{-7/4}}{8}\left(-\frac{2\pi^3}{15}-\frac{\pi z^2\beta}{3}+\frac{z^4\beta^2}{20\pi}\right),
%\end{align}
%and so on. 
Now let $j$ be even, say $j=2n$. Then using the fact \cite[p.~5, Equation (1.14)]{temme} that
\begin{equation*}
\frac{(-1)^{m-1}2^{2m-1}\pi^{2m}}{(2m)!}B_{2m}=\zeta(2m)
\end{equation*}
in the second step below, we see that
{\allowdisplaybreaks\begin{align}\label{evenc}
\sum_{n=0}^{\infty}I_{2n,z}&=\frac{2}{z}\beta^{\frac{1}{4}}e^{\frac{z^2}{8}}\sum_{n=0}^{\infty}\frac{\Gamma\left(n+\frac{1}{2}\right)}{\beta^{n+\frac{1}{2}}}\sum_{m=0}^{n}\frac{B_{2n-2m}(2\pi)^{2n-2m-1}}{(2n-2m)!}\frac{(-1)^{m}(z\sqrt{\beta})^{2m+1}}{(2m+1)!}\nonumber\\
&=\frac{2}{z}\beta^{\frac{1}{4}}e^{\frac{z^2}{8}}\sum_{n=0}^{\infty}\frac{(-1)^n(z\sqrt{\b})^{2n+1}}{\beta^{n+\frac{1}{2}}}\Gamma\left(n+\frac{1}{2}\right)\sum_{m=0}^{n}\frac{B_{2m}(2\pi)^{2m-1}}{(2m)!}\frac{(-1)^{m}(z\sqrt{\beta})^{-2m}}{(2n-2m+1)!}\nonumber\\
&=\frac{-2}{\pi z}\beta^{\frac{1}{4}}e^{\frac{z^2}{8}}\sum_{n=0}^{\infty}(-1)^nz^{2n+1}\Gamma\left(n+\frac{1}{2}\right)\sum_{m=0}^{n}\frac{\zeta(2m)}{(z\sqrt{\beta})^{2m}(2n-2m+1)!}\nonumber\\
&=\frac{-2}{\pi}\beta^{\frac{1}{4}}e^{\frac{z^2}{8}}\sum_{m=0}^{\infty}\frac{\zeta(2m)}{(z\sqrt{\beta})^{2m}}\sum_{n=m}^{\infty}\frac{(-1)^nz^{2n}}{(2n-2m+1)!}\G\left(n+\frac{1}{2}\right)\nonumber\\
&=\frac{-2}{\pi}\beta^{\frac{1}{4}}e^{\frac{z^2}{8}}\sum_{m=0}^{\infty}\frac{(-1)^m\zeta(2m)}{\b^{m}}\G\left(m+\frac{1}{2}\right){}_1F_{1}\left(m+\frac{1}{2};\frac{3}{2};\frac{-z^2}{4}\right).\nonumber\\
%&=\frac{-2}{\sqrt{\pi}}\a^{-1/4}e^{z^2/8}\sum_{m=0}^{\infty}\left(\frac{-\a}{\pi^2}\right)^{m}\zeta(2m)\G\left(m+\frac{1}{2}\right){}_1F_{1}\left(m+\frac{1}{2};\frac{3}{2};\frac{-z^2}{4}\right),
\end{align}}
From \eqref{ibetaz}, \eqref{a14} and \eqref{evenc}, we obtain the asymptotic expansion of $I(iz,\b)$ as $\b\to\infty$ as
\begin{align}\label{ibetaz1}
I(iz, \beta)&\sim\left(\frac{\sqrt{\pi}}{z}e^{-\frac{z^2}{8}}\textup{erfi}\left(\frac{z}{2}\right)-e^{\frac{z^2}{8}}\sum_{m=0}^{\infty}\frac{(-1)^mm!z^{2m}}{(2m+1)!}\right)\beta^{-\frac{1}{4}}\nonumber\\
&\quad-\frac{2}{\pi}\beta^{\frac{1}{4}}e^{\frac{z^2}{8}}\sum_{m=0}^{\infty}\frac{(-1)^m\zeta(2m)}{\b^{m}}\G\left(m+\frac{1}{2}\right){}_1F_{1}\left(m+\frac{1}{2};\frac{3}{2};\frac{-z^2}{4}\right).
\end{align}
%\begin{align}\label{alphaasym0}
%&\left(\frac{ze^{z^2/8}}{\sqrt{\pi}}-\frac{z^3e^{z^2/8}}{12\sqrt{\pi}}+\frac{z^5e^{z^2/8}}{160\sqrt{\pi}}-\frac{z^7e^{z^2/8}}{2688\sqrt{\pi}}+\frac{z^9e^{z^2/8}}{55296\sqrt{\pi}}-\cdots\right)\alpha^{-1/4}\nonumber\\
%&\left\{\frac{e^{-z^2/8}}{z}\textup{erfi}\left(\frac{z}{2}\right)-\frac{e^{z^2/8}}{z\sqrt{\pi}}\sum_{n=0}^{\infty}\frac{(-1)^nz^{2n+1}n!}{(2n+1)!}\right\}\alpha^{1/4}+\frac{2e^{z^2/8}}{z}\left(\sum_{n=0}^{\infty}\frac{(-1)^n(z/2)^{2n+1}}{(2n+1)n!}\right)\alpha^{-1/4}\nonumber\\
%-\frac{ze^{z^2/8}}{\pi}+\frac{z^3e^{z^2/8}}{6\pi}\right\}\alpha^{1/4}\nonumber\\
%&\quad+\frac{e^{z^2/8}}{6}\left(\sum_{n=0}^{\infty}\frac{(-1)^nz^{2n}}{4^nn!}\right)\alpha^{3/4}-\frac{e^{z^2/8}}{45\sqrt{\pi}}\left(\sum_{n=0}^{\infty}\frac{\Gamma\left(n+\frac{5}{2}\right)(-z^2)^n}{(2n+1)!}\right)\alpha^{7/4}\nonumber\\
%&\quad+\frac{2e^{z^2/8}}{945\sqrt{\pi}}\left(\sum_{n=0}^{\infty}\frac{\Gamma\left(n+\frac{7}{2}\right)(-z^2)^n}{(2n+1)!}\right)\alpha^{11/4}-\cdots.
%\end{align}
Note that $\textup{erf}(z)$ has the following Taylor series expansion, which is valid for all $z\in\mathbb{C}$ \cite[p.~162, 7.6.1]{nist}:
\begin{equation}\label{erfseries}
\textup{erf}(z)=\frac{2}{\sqrt{\pi}}\sum_{n=0}^{\infty}\frac{(-1)^nz^{2n+1}}{n!(2n+1)}.
\end{equation}
From \eqref{errrel},
\begin{equation}\label{erfiseries}
\textup{erfi}(z)=\frac{2}{\sqrt{\pi}}\sum_{n=0}^{\infty}\frac{z^{2n+1}}{n!(2n+1)}.
\end{equation}
It is now easy to see that
\begin{equation}\label{giv0}
\sum_{m=0}^{\infty}\frac{(-1)^mm!z^{2m}}{(2m+1)!}=\frac{\sqrt{\pi}}{z}e^{-z^2/4}\textup{erfi}\left(\frac{z}{2}\right).
\end{equation}
%Also all of the sums in \eqref{alphaasym0} have closed-form expressions that can be easily evaluated.
%\begin{equation}\label{giv1}
%\frac{z}{6}-\frac{z^3}{24}+\frac{z^5}{192}-\frac{z^7}{2304}+\cdots=\frac{z}{6}e^{-z^2/4}.
%\end{equation}
Substituting \eqref{giv0} in \eqref{ibetaz1}, we arrive at 
\begin{align*}
I(iz, \beta)&\sim-\frac{2}{\pi}\beta^{\frac{1}{4}}e^{\frac{z^2}{8}}\sum_{m=0}^{\infty}\frac{(-1)^m\zeta(2m)}{\b^{m}}\G\left(m+\frac{1}{2}\right){}_1F_{1}\left(m+\frac{1}{2};\frac{3}{2};\frac{-z^2}{4}\right).
\end{align*}
%This completes the proof of Theorem \ref{1asy}.
Since $\a\b=\pi^2$, using \eqref{gen}, this also gives the asymptotic expansion of $I(z,\a)$ as $\a\to 0$ as claimed in \eqref{1asyeq}.
\end{proof}

The second error function transformation in Theorem \ref{thmseft} can be rephrased for $\a\b=\pi^2$ and $z\neq 0$ as follows:
\begin{align}\label{sefty}
K(z,\a)&:=\frac{\sqrt{\pi}}{z}\alpha^{-\frac{1}{4}}e^{\frac{z^2}{8}}\textup{erf}\left(\frac{z}{2}\right)-\frac{4}{z}\alpha^{\frac{1}{4}}e^{-\frac{z^2}{8}}\int_{-\infty}^{0}\frac{e^{-\alpha x^2}\sinh(\sqrt{\a} xz)}{e^{2\pi x}-1}\, dx\nonumber\\
&=\frac{\sqrt{\pi}}{z}\b^{-\frac{1}{4}}e^{-\frac{z^2}{8}}\textup{erfi}\left(\frac{z}{2}\right)-\frac{4}{z}\b^{\frac{1}{4}}e^{\frac{z^2}{8}}\int_{-\infty}^{0}\frac{e^{-\b x^2}\sin(\sqrt{\b} xz)}{e^{2\pi x}-1}\, dx=:K(iz,\b).
\end{align}
The analogue of Theorem \ref{1asy} for the above identity is given below. The details are similar to those in the proof of Theorem \ref{1asy} and hence not provided.
\begin{theorem}\label{1asy2}
Fix $z\in\mathbb{C}$. As $\a\to 0$,
\begin{align*}
K(z,\a)&\sim\frac{2}{\sqrt{\pi}}\a^{-1/4}e^{z^2/8}\sum_{m=0}^{\infty}\left(\frac{-\a}{\pi^2}\right)^{m}\zeta(2m)\G\left(m+\frac{1}{2}\right){}_1F_{1}\left(m+\frac{1}{2};\frac{3}{2};\frac{-z^2}{4}\right).
\end{align*}
%\begin{align*}
%K(z,\a)&\sim-\frac{\sqrt{\pi}}{z}e^{z^2/8}\textup{erf}\left(\frac{z}{2}\right)\alpha^{-1/4}-\frac{e^{-z^2/8}}{6}\alpha^{3/4}-\frac{(z^2-6)e^{-z^2/8}}{360}\alpha^{7/4}\nonumber\\
%&\quad-\frac{(60-20z^2+z^4)e^{-z^2/8}}{15120}\alpha^{11/4}-\cdots.
%\end{align*}
\end{theorem}
Combining Theorems \ref{1asy} and \ref{1asy2} leads us to the following asymptotic expansion of the integral analogue of theta function.
\begin{theorem}\label{ysa}
Fix $z\in\mathbb{C}$. As $\a\to 0$,
\begin{align*}
&\int_{-\infty}^{\infty}\frac{e^{-\a x^2}\sinh(\sqrt{\a}xz)}{e^{2\pi x}-1}\, dx\nonumber\\
&\sim\frac{-z}{\sqrt{\pi\a}}e^{z^2/4}\sum_{m=0}^{\infty}\left(\frac{-\a}{\pi^2}\right)^{m}\zeta(2m)\G\left(m+\frac{1}{2}\right){}_1F_{1}\left(m+\frac{1}{2};\frac{3}{2};\frac{-z^2}{4}\right).
\end{align*}
\end{theorem}
\section{Generalization of Ramanujan's approximation and integral identities involving hypergeometric functions}
The results from this section follow from successively differentiating the first error function transformation in the form \eqref{gen} with respect to $z$. The presence of $-x^2$ in the exponential term in the numerator of either sides justifies differentiation under the integral sign. As mentioned in the introduction, differentiating \eqref{gen} $n$ times where $n$ is odd just gives $0=0$. However when $n$ is even, two different behaviors are observed accordingly as $n\equiv0\hspace{1mm}(\text{mod}\hspace{1mm} 4)$ and as $n\equiv2\hspace{1mm}(\text{mod}\hspace{1mm} 4)$.

\subsection{The case $n\equiv0\hspace{1mm}(\text{mod}\hspace{1mm} 4)$}
\begin{proof}[Theorem \textup{\ref{thmtrn0m4}}][]
Let $n=4k$ for a non-negative integer $k$. Let 
\begin{equation}\label{ga}
G(z,\a):=\frac{\sqrt{\pi}}{z}\alpha^{-\frac{1}{4}}e^{\frac{z^2}{8}}\textup{erf}\left(\frac{z}{2}\right)+\frac{4}{z}\alpha^{\frac{1}{4}}e^{-\frac{z^2}{8}}\int_{0}^{\infty}\frac{e^{-\alpha x^2}\sinh(\sqrt{\alpha} xz)}{e^{2\pi x}-1}\, dx
\end{equation}
and let
\begin{equation}\label{ika}
H_{k}(\a):=\left.\frac{d^{4k}}{dz^{4k}}G(z,\a)\right|_{z=0}.
\end{equation}
Using \eqref{gen} leads us to a transformation of the form $H_{k}(\a)=H_{k}(\b)$ since 
\begin{equation}\label{ikab}
H_{k}(\a)=\left.\frac{d^{4k}}{dz^{4k}}G(z,\a)\right|_{z=0}=\left.\frac{d^{4k}}{d(iz)^{4k}}G(iz,\b)\right|_{z=0}=H_{k}(\b).
\end{equation}
Multiply the power series expansion of $e^{z^2/8}$ with that of $\textup{erf}(z)$ given in \eqref{erfseries} and extract the coefficient of $z^{4k}$ in the product. This gives
\begin{align*}
\left.\frac{d^{4k}}{dz^{4k}}\frac{e^{\frac{z^2}{8}}}{z}\textup{erf}\left(\frac{z}{2}\right)\right|_{z=0}&=\frac{2(4k)!}{\sqrt{\pi}}\sum_{m=0}^{2k}\frac{(-1)^m}{2^{4k+2m+1}m!(2k-m)!(4k-2m+1)}\nonumber\\
&=\frac{(4k)!}{2^{4k}\sqrt{\pi}(2k)!(4k+1)}{}_2F_{1}\left(-\frac{1}{2}-2k,-2k;\frac{1}{2}-2k;\frac{1}{2}\right).
\end{align*}
In \cite[p.~113, Equation (5.12)]{temme}, we find the following hypergeometric transformation, valid for $|\arg(1-z)|<\pi$:
\begin{align}\label{transformation}
{}_2F_{1}(a,b;c;z)&=\frac{\G(c)\G(b-a)}{\G(b)\G(c-a)}(1-z)^{-a}{}_2F_{1}\left(a,c-b;a-b+1;\frac{1}{1-z}\right)\nonumber\\
&\quad+\frac{\G(c)\G(a-b)}{\G(a)\G(c-b)}(1-z)^{-b}{}_2F_{1}\left(b,c-a;b-a+1;\frac{1}{1-z}\right).
\end{align}
Use this transformation with $z=1/2, a=-\frac{1}{2}-2k, b=-2k, c=\frac{1}{2}-2k$ to obtain
\begin{equation*}
{}_2F_{1}\left(-\frac{1}{2}-2k,-2k;\frac{1}{2}-2k;\frac{1}{2}\right)=\frac{(4k+1)}{2^{2k}}{}_2F_{1}\left(-2k,1;\frac{3}{2};2\right).
\end{equation*}
This gives
\begin{align}\label{ika1}
\left.\frac{d^{4k}}{dz^{4k}}\frac{e^{\frac{z^2}{8}}}{z}\textup{erf}\left(\frac{z}{2}\right)\right|_{z=0}=\frac{(4k)!}{2^{6k}\sqrt{\pi}(2k)!}{}_2F_{1}\left(-2k,1;\frac{3}{2};2\right).
\end{align}
Similarly multiplying the power series expansion of $e^{-z^2/8}$ with that of $\sinh(\sqrt{\a}xz)$ and extracting the coefficient of $z^{4k}$ in the product, we get
\begin{align*}
\left.\frac{d^{4k}}{dz^{4k}}\frac{e^{-\frac{z^2}{8}}}{z}\sinh(\sqrt{\a}xz)\right|_{z=0}&=(4k)!\sum_{m=0}^{2k}\frac{(-1)^m(\sqrt{\a}x)^{4k-2m+1}}{8^mm!(4k-2m+1)!}\nonumber\\
&=\frac{(4k)!\sqrt{\a}x}{2^{6k}}\sum_{m=0}^{2k}\frac{(-8)^{m}(\sqrt{\a}x)^{2m}}{(2k-m)!(2m+1)!},
\end{align*}
where we replaced $m$ by $2k-m$ in the last step. It can be seen that
\begin{equation*}
\sum_{m=0}^{2k}\frac{(-8)^{m}(\sqrt{\a}x)^{2m}}{(2k-m)!(2m+1)!}=\frac{1}{(2k)!}{}_1F_{1}\left(-2k;\frac{3}{2};2\a x^2\right).
\end{equation*}
Thus
\begin{align}\label{ika2}
\left.\frac{d^{4k}}{dz^{4k}}\frac{e^{-\frac{z^2}{8}}}{z}\sinh(\sqrt{\a}xz)\right|_{z=0}=\frac{(4k)!\sqrt{\a}x}{2^{6k}(2k)!}{}_1F_{1}\left(-2k;\frac{3}{2};2\a x^2\right).
\end{align}
Hence from \eqref{ga}, \eqref{ika}, \eqref{ika1} and \eqref{ika2}, we obtain
\begin{equation*}
H_{k}(\a)=\frac{(4k)!}{2^{6k}(2k)!}\left\{\a^{-1/4}{}_2F_{1}\left(-2k,1;\frac{3}{2};2\right)+4\a^{3/4}\int_{0}^{\infty}\frac{xe^{-\a x^2}}{e^{2\pi x}-1}{}_1F_{1}\left(-2k;\frac{3}{2};2\a x^2\right)\, dx\right\},
\end{equation*}
which when combined with \eqref{ikab}, proves Theorem \ref{thmtrn0m4}.
\end{proof}
Next, we prove a generalization of Ramanujan's approximation in \eqref{near}.
\begin{proof}[Theorem \textup{\ref{thmgennear}}][]
Let $H_{k}(\b)$ be as defined in \eqref{ika} and consider the integral 
\begin{equation}
\int_{0}^{\infty}\frac{xe^{-\b x^2}}{e^{2\pi x}-1}{}_1F_{1}\left(-2k;\frac{3}{2};2\b x^2\right)\, dx.
\end{equation}
Employing the change of variable $x=\sqrt{t}$, using the series definition of ${}_1F_{1}$ and interchanging the order of summation and integration, it is seen that
\begin{equation*}
\int_{0}^{\infty}\frac{xe^{-\b x^2}}{e^{2\pi x}-1}{}_1F_{1}\left(-2k;\frac{3}{2};2\b x^2\right)\, dx=\frac{1}{2}\sum_{m=0}^{2k}\frac{(-2k)_{m}(2\b)^{m}}{\left(\frac{3}{2}\right)_{m}m!}\int_{0}^{\infty}\frac{e^{-\b t}t^m}{e^{2\pi\sqrt{t}}-1}\, dt.
\end{equation*}
Now use the generating function for Bernoulli numbers to obtain, for $|t|<1$,
\begin{equation*}
\frac{t^m}{e^{2\pi\sqrt{t}}-1}=\sum_{j=0}^{\infty}\frac{B_{j}(2\pi)^{j-1}}{j!}t^{\frac{j-1+2m}{2}}.
\end{equation*}
Using Watson's lemma from \eqref{watlem1}, we find that as $\b\to\infty$,
\begin{align*}
H_{k}(\b)\sim\b^{-1/4}{}_2F_{1}\left(-2k,1;\frac{3}{2};2\right)+2\b^{3/4}\sum_{m=0}^{2k}\frac{(-2k)_{m}(2\b)^{m}}{\left(\frac{3}{2}\right)_{m}m!}\sum_{j=0}^{\infty}\frac{B_{j}(2\pi)^{j-1}}{j!\b^{\frac{j+2m+1}{2}}}\G\left(\frac{j+2m+1}{2}\right).
\end{align*}
Since $H_{k}(\a)=H_{k}(\b)$, using the fact $\b=\pi^2/\a$ yields for $\a\to 0$,
\begin{align*}
H_{k}(\a)&\sim\frac{\a^{1/4}}{\sqrt{\pi}}{}_2F_{1}\left(-2k,1;\frac{3}{2};2\right)+\frac{2}{\sqrt{\pi}}\a^{-1/4}\sum_{m=0}^{2k}\frac{(-2k)_{m}2^{m}}{\left(\frac{3}{2}\right)_{m}m!}\sum_{j=0}^{\infty}\frac{B_{j}2^{j-1}\a^{j/2}}{j!}\G\left(\frac{j+2m+1}{2}\right)\nonumber\\
&=\frac{\a^{-1/4}}{\sqrt{\pi}}\sum_{j=0\atop{j\neq 1}}^{\infty}\frac{B_{j}2^{j}\a^{j/2}}{j!}\G\left(\frac{j+1}{2}\right){}_2F_{1}\left(-2k,\frac{1+j}{2};\frac{3}{2};2\right)\nonumber\\
&=\a^{-1/4}{}_2F_{1}\left(\frac{1}{2},-2k;\frac{3}{2};2\right)+\frac{\a^{3/4}}{6}+\frac{\a^{-1/4}}{\sqrt{\pi}}\sum_{j=3}^{\infty}\frac{B_{j}2^{j}\a^{j/2}}{j!}\G\left(\frac{j+1}{2}\right){}_2F_{1}\left(-2k,\frac{1+j}{2};\frac{3}{2};2\right).
\end{align*}
We now find a simpler function, namely the one claimed on the right-hand side of \eqref{gennear}, that is ``nearly'' equal to $H_{k}(\a)$ when $\a$ is very small in the sense that the asymptotic expansion of this simpler function as $\a\to 0$ matches the first two terms in those of $H_{k}(\a)$. Note that such a function should preserve the invariance under replacing $\a$ by $\b$ and vice-versa. In order to match the leading term in the asymptotic expansion, we raise $1/\a$ to the power $1/4$ and have its coefficient as $\displaystyle {}_2F_{1}\left(\tfrac{1}{2},-2k;\tfrac{3}{2};2\right)$, which is equal to ${}_2F_{1}\left(-2k,1;\frac{3}{2};2\right)$ by Pfaff's transformation \cite[p.~110, Equation (5.5)]{temme}
\begin{equation*}\label{pfaff}
{}_2F_{1}(a,b;c;z)=(1-z)^{-b}{}_2F_{1}\left(c-a,b;c;\frac{z}{z-1}\right).
\end{equation*}
Since the approximating function has to be symmetric, we need to raise $1/\b$ along with $1/\a$ to the power $1/4$. So the function we are seeking assumes the form
\begin{equation*}
{}_2F_{1}\left(-2k,1;\frac{3}{2};2\right)\left(\frac{1}{\a}+\frac{1}{\b}+c(k)\right)^{1/4},
\end{equation*}
where $c(k)$ is some constant depending on only $k$. Since $\a$ is very small, the main contribution in the asymptotic expansion comes from $1/\a$ and $c(k)$ but not from $1/\b$. Thus, the next term in the Taylor series of this function is
\begin{equation*}
\frac{\a^{3/4}}{4}{}_2F_{1}\left(-2k,1;\frac{3}{2};2\right)c(k).
\end{equation*}
As we want this to be equal to $\frac{\a^{3/4}}{6}$, it is clear that $c(k)=\frac{2}{3\cdot{}_2F_{1}\left(-2k,1;\frac{3}{2};2\right)}$.  Thus the required function is
\begin{equation*}
{}_2F_{1}\left(-2k,1;\frac{3}{2};2\right)\left(\frac{1}{\a}+\frac{1}{\b}+\frac{2}{3\cdot{}_2F_{1}\left(-2k,1;\frac{3}{2};2\right)}\right)^{1/4}.
\end{equation*}
This completes the proof of \eqref{gennear}.
\end{proof}

\subsection{The case $n\equiv2\hspace{1mm}(\text{mod}\hspace{1mm} 4)$}
\begin{proof}[Theorem \textup{\ref{thmtrn2m4}}][]
Now let $n=4k+2$, where $k$ is again any non-negative integer. Let
\begin{equation*}
J_{k}(\a):=\left.\frac{d^{4k+2}}{dz^{4k+2}}G(z,\a)\right|_{z=0},
\end{equation*}
where $G(z,\a)$ is defined in \eqref{ga}. This time \eqref{gen} gives us a transformation of the form $J_{k}(\a)=-J_{k}(\b)$ since 
\begin{equation*}
J_{k}(\a)=\left.\frac{d^{4k+2}}{dz^{4k+2}}G(z,\a)\right|_{z=0}=\left.\frac{d^{4k+2}}{dz^{4k+2}}G(iz,\b)\right|_{z=0}=-\left.\frac{d^{4k+2}}{d(iz)^{4k+2}}G(iz,\b)\right|_{z=0}=-J_{k}(\b).
\end{equation*}
The proof is now similar to that of Theorem \ref{thmtrn0m4} and so we state only the important steps. Let us start with the fact that 
\begin{equation}\label{jka0}
-J_{k}(\b)=\left.\frac{d^{4k+2}}{dz^{4k+2}}G(iz,\b)\right|_{z=0}.
\end{equation}
Multiply the power series expansion of $e^{-z^2/8}$ with that of $\textup{erfi}\left(\frac{z}{2}\right)$ given in \eqref{erfiseries} and extract the coefficient of $z^{4k+2}$ in the product. Then identifying the coefficient as a hypergeometric function and using the transformation \eqref{transformation} with $z=1/2, a=-\frac{3}{2}-2k, b=-1-2k, c=-\frac{1}{2}-2k$ to simplify this hypergeometric function results in
\begin{align*}
\left.\frac{d^{4k+2}}{dz^{4k+2}}\frac{e^{\frac{-z^2}{8}}}{z}\textup{erfi}\left(\frac{z}{2}\right)\right|_{z=0}=\frac{(4k+2)!}{2^{6k+3}\sqrt{\pi}(2k+1)!}{}_2F_{1}\left(-2k-1,1;\frac{3}{2};2\right).
\end{align*}
Similarly,
{\allowdisplaybreaks\begin{align}\label{jka2}
\left.\frac{d^{4k+2}}{dz^{4k+2}}\frac{e^{\frac{z^2}{8}}}{z}\sin(\sqrt{\b}xz)\right|_{z=0}&=-(4k+2)!\sum_{m=0}^{2k+1}\frac{(-1)^m(\sqrt{\b}x)^{4k-2m+3}}{8^mm!(4k-2m+3)!}\nonumber\\
&=\frac{(4k+2)!\sqrt{\b}x}{2^{6k+3}}\sum_{m=0}^{2k+1}\frac{(-8)^{m}(\sqrt{\b}x)^{2m}}{(2k+1-m)!(2m+1)!}\nonumber\\
&=\frac{(4k+2)!\sqrt{\b}x}{2^{6k+3}(2k+1)!}{}_1F_{1}\left(-2k-1;\frac{3}{2};2\b x^2\right).
\end{align}}
From \eqref{jka0} and \eqref{jka2}, 
\begin{align*}
-J_{k}(\b)&=\frac{(4k+2)!}{2^{6k+3}(2k+1)!}\bigg\{\b^{-1/4}{}_2F_{1}\left(-2k-1,1;\frac{3}{2};2\right)\nonumber\\
&\quad\quad\quad\quad\quad\quad\quad\quad+4\b^{3/4}\int_{0}^{\infty}\frac{xe^{-\b x^2}}{e^{2\pi x}-1}{}_1F_{1}\left(-2k-1;\frac{3}{2};2\b x^2\right)\, dx\bigg\}.
\end{align*}
This proves Theorem \ref{thmtrn2m4}.
\end{proof}
As remarked in the introduction, results corresponding to the ones in Theorem \ref{thmtrn0m4} - Corollary \ref{corexact} can be obtained using similar techniques through \eqref{sefty}. These results are collectively put in the following theorem. We refrain from giving the proof since the details are similar to those of Theorem \ref{thmtrn0m4} - Corollary \ref{corexact}.
\begin{theorem}\label{appseft}
Let $\a, \b$ be two positive numbers such that $\a\b=\pi^2$. Then
{\allowdisplaybreaks\begin{align}
&\textup{(i)}\hspace{2mm} \a^{-1/4}{}_2F_{1}\left(-2k,1;\frac{3}{2};2\right)-4\a^{3/4}\int_{-\infty}^{0}\frac{xe^{-\a x^2}}{e^{2\pi x}-1}{}_1F_{1}\left(-2k;\frac{3}{2};2\a x^2\right)\, dx\nonumber\\
&\quad\quad=\b^{-1/4}{}_2F_{1}\left(-2k,1;\frac{3}{2};2\right)-4\b^{3/4}\int_{-\infty}^{0}\frac{xe^{-\b x^2}}{e^{2\pi x}-1}{}_1F_{1}\left(-2k;\frac{3}{2};2\b x^2\right)\, dx.\nonumber\\
&\textup{(ii)}\hspace{2mm} \a^{-1/4}{}_2F_{1}\left(-2k,1;\frac{3}{2};2\right)-4\a^{3/4}\int_{-\infty}^{0}\frac{xe^{-\a x^2}}{e^{2\pi x}-1}{}_1F_{1}\left(-2k;\frac{3}{2};2\a x^2\right)\, dx\nonumber\\
&\quad\quad=-{}_2F_{1}\left(-2k,1;\frac{3}{2};2\right)\left(\frac{1}{\a}+\frac{1}{\b}+\frac{2}{3\cdot{}_2F_{1}\left(-2k,1;\frac{3}{2};2\right)}\right)^{1/4}, \text{``nearly''}.\nonumber\\
&\textup{(iii)}\hspace{2mm} -\a^{-1/4}{}_2F_{1}\left(-2k-1,1;\frac{3}{2};2\right)+4\a^{3/4}\int_{-\infty}^{0}\frac{xe^{-\a x^2}}{e^{2\pi x}-1}{}_1F_{1}\left(-2k-1;\frac{3}{2};2\a x^2\right)\, dx\nonumber\\
&=\b^{-1/4}{}_2F_{1}\left(-2k-1,1;\frac{3}{2};2\right)-4\b^{3/4}\int_{-\infty}^{0}\frac{xe^{-\b x^2}}{e^{2\pi x}-1}{}_1F_{1}\left(-2k-1;\frac{3}{2};2\b x^2\right)\, dx.\nonumber
\end{align}}
In particular when $\a=\b=\pi$, we have
\begin{equation}\label{exact1}
\int_{-\infty}^{0}\frac{xe^{-\pi x^2}}{e^{2\pi x}-1}{}_1F_{1}\left(-2k-1;\frac{3}{2};2\pi x^2\right)\, dx=\frac{1}{4\pi}\cdot{}_2F_{1}\left(-2k-1,1;\frac{3}{2};2\right).
\end{equation}
\end{theorem}
\section{Generalization of an asymptotic expansion of Oloa}
We first explain how the integral in \eqref{ramanujan-98} is related to the one on the extreme right side of \eqref{mrram}. Write the latter integral as
\begin{equation}\label{amrram}
\int_{0}^{\infty}(1+t^2)\Gamma\left(\frac{-1+it}{4}\right)\Gamma\left(\frac{-1-it}{4}\right)\frac{\Xi\left(\frac{t}{2}\right)}{1+t^2}\cos \left(\frac{1}{2}t\log\alpha\right)\, dt.
\end{equation}
If we now square the expression $\frac{\Xi\left(\frac{t}{2}\right)}{1+t^2}$ in \eqref{amrram}, then as discussed in \cite{dixitmoll}, this amounts to squaring the functional equation of the Riemann zeta function, and moreover the squared expression admits a generalization 
\begin{equation*}
\frac{ 
\Xi \left( \frac{t + i z}{2} \right) 
\Xi \left( \frac{t - i z}{2} \right) }
{ (t^{2} + (z+1)^{2})(t^{2}+ (z-1)^{2})}.
\end{equation*}
This is what Ramanujan may have had at the back of his mind when he worked \cite[Section 5]{riemann} with the generalization 
\begin{equation*}
\int_{0}^{\infty}(t^2+(z-1)^2)\Gamma\left(\frac{z-1+it}{4}\right)\Gamma\left(\frac{z-1-it}{4}\right)\frac{\Xi \left( \frac{t + i z}{2} \right) 
\Xi \left( \frac{t - i z}{2} \right) }
{ (t^{2} + (z+1)^{2})(t^{2}+ (z-1)^{2})}\cos \left(\frac{1}{2}t\log\alpha\right)\, dt,
\end{equation*}
of \eqref{amrram}, which upon simplification gives \eqref{ramanujan-98}. Ramanujan \cite{riemann} obtains a transformation formula associated with this integral. In \cite{dixitmoll}, Moll and one of the authors found the following new representation of this transformation, which generalizes a transformation of Koshliakov \cite[Equation (6)]{koshli} \cite[Equations (21), (27)]{kosh1937}.\\ 

\textit{Assume $-1 <$ \textup{Re} $z < 1$. Let $\Omega(x, z)$ be defined by
\begin{equation*}
\Omega(x,z) = 2 \sum_{n=1}^{\infty} \sigma_{-z}(n) n^{z/2} 
\left( e^{\pi i z/4} K_{z}( 4 \pi e^{\pi i/4} \sqrt{nx} ) +
 e^{-\pi i z/4} K_{z}( 4 \pi e^{-\pi i/4} \sqrt{nx} ) \right),
\end{equation*}
where $\sigma_{-z}(n)=\sum_{d|n}d^{-z}$ and $K_\nu(z)$ denotes the modified Bessel function of order $\nu$. Then for $\alpha, \beta>0, \alpha\beta=1$, 
\begin{align}\label{genelkosh}
&\alpha^{(z+1)/2} 
\int_{0}^{\infty} e^{-2 \pi \alpha x} x^{z/2} 
\left( \Omega(x,z) - \frac{1}{2 \pi} \zeta(z) x^{z/2-1} \right) dx\nonumber \\
&=\beta^{(z+1)/2} 
\int_{0}^{\infty} e^{-2 \pi \beta x} x^{z/2} 
\left( \Omega(x,z) - \frac{1}{2 \pi} \zeta(z) x^{z/2-1} \right) dx\nonumber\\
&=\frac{1}{2\pi^{(z+5)/2}}
\int_{0}^{\infty} 
\Gamma \left( \frac{z-1+it}{4} \right)
\Gamma \left( \frac{z-1-it}{4} \right)
\Xi \left( \frac{t + iz}{2} \right)
\Xi \left( \frac{t - iz}{2} \right)\frac{\cos( \tfrac{1}{2} t \log \alpha) \, dt}{(t^{2} + (z+1)^{2})}.
\end{align}}

Theorem \ref{2asy} is now proved using \eqref{genelkosh}.
\begin{proof}[Theorem \textup{\ref{2asy}}][]
We obtain the asymptotic expansion of the integral indirectly by obtaining the same for the integral on the extreme left of \eqref{genelkosh}. Let 
\begin{equation*}
g(t, z):=t^{z/2} 
\left( \Omega(t,z) - \frac{1}{2 \pi} \zeta(z) t^{z/2-1} \right).
\end{equation*}
%We now apply Watson's lemma to the integral 
%\begin{equation}
%\int_{0}^{\infty}e^{-2\pi\a t}
%\end{equation}
%from \eqref{genelkosh}. 
%As per Watson's lemma, we need to find $a_m$ and constants $\lambda$ and $\mu$ with Re $\lambda>0$ and $\mu>0$ so that as $t\to 0$,
%\begin{equation}
%g(t, z)\sim\sum_{m=0}^{\infty}a_mt^{(m+\lambda-\mu)/\mu}.
%\end{equation}
We use the following identity  established in \cite[Proposition 6.1]{dixitmoll}.
\begin{equation*}
\Omega(t,z) = - \frac{\Gamma(z) \zeta(z)}{(2 \pi \sqrt{t})^{z}}  +
\frac{t^{z/2-1}}{2 \pi} \zeta(z) - 
\frac{t^{z/2}}{2} \zeta(z+1)  +
\frac{t^{z/2+1}}{\pi} \sum_{n=1}^{\infty} \frac{\sigma_{-z}(n)}{n^{2}+t^{2}}.
\end{equation*}
Since for $|t|<1$,
{\allowdisplaybreaks\begin{align*}
\sum_{n=1}^{\infty} \frac{\sigma_{-z}(n)}{n^{2}+t^{2}}&=\sum_{n=1}^{\infty} \frac{\sigma_{-z}(n)}{n^2}\frac{1}{\left(1+\left(\frac{t}{n}\right)^{2}\right)}\nonumber\\
&=\sum_{n=1}^{\infty}\frac{\sigma_{-z}(n)}{n^2}\sum_{m=0}^{\infty}(-1)^m\left(\frac{t}{n}\right)^{2m}\nonumber\\
&=\sum_{m=0}^{\infty}(-1)^mt^{2m}\sum_{n=1}^{\infty}\frac{\sigma_{-z}(n)}{n^{2m+2}}\nonumber\\
&=\sum_{m=0}^{\infty}(-1)^m\zeta(2m+2)\zeta(2m+2+z)t^{2m},
\end{align*}}%
we see that if
\begin{equation*}
h(t,z):=g(t,z)+\frac{\Gamma(z) \zeta(z)}{(2 \pi)^{z}}+\frac{t^{z}}{2} \zeta(z+1),
\end{equation*}
then
\begin{equation*}
h(t,z)=\sum_{m=0}^{\infty}\frac{(-1)^m}{\pi}\zeta(2m+2)\zeta(2m+2+z)t^{2m+z+1},
\end{equation*}
so also as $t\to 0$,
\begin{equation*}
h(t,z)\sim\sum_{m=0}^{\infty}\frac{(-1)^m}{\pi}\zeta(2m+2)\zeta(2m+2+z)t^{2m+z+1}.
\end{equation*}

We now apply Lemma \ref{watlem} with $\lambda=(z+2)/2$ and $\mu=1/2$. The condition $-1<$ Re $z<1$ guarantees that Re $\lambda>0$ which is necessary as remarked after Lemma \ref{watlem}. Then as $\a\to\infty$,
\begin{equation}\label{aftwat1}
\int_{0}^{\infty}e^{-2\pi\a t}h(t, z)\, dt\sim\sum_{m=0}^{\infty}\frac{(-1)^m}{\pi(2\pi\a)^{2m+z+2}}\G(2m+2+z)\zeta(2m+2)\zeta(2m+2+z).
\end{equation}
Note that
\begin{align}\label{aftwat2}
\int_{0}^{\infty}e^{-2\pi\a t}h(t, z)\, dt&=\int_{0}^{\infty}e^{-2\pi\a t}g(t, z)\, dt+\frac{\Gamma(z) \zeta(z)}{(2 \pi)^{z}}\int_{0}^{\infty}e^{-2\pi\a t}\, dt\nonumber\\
&\quad+\frac{\zeta(z+1)}{2}\int_{0}^{\infty}e^{-2\pi\a t}t^{z}\, dt\nonumber\\
&=\int_{0}^{\infty}e^{-2\pi\a t}g(t, z)\, dt+\frac{\Gamma(z) \zeta(z)}{\a(2 \pi )^{z+1}}+\frac{\zeta(z+1)\G(z+1)}{2(2 \pi\a )^{z+1}}.
\end{align}
From \eqref{aftwat1} and \eqref{aftwat2} and \eqref{genelkosh}, one now obtains \eqref{oloag} after simplification.
\end{proof}

\section{Concluding remarks and some open questions}
1. In this paper, we found two new transformations involving error functions, namely the ones in Theorems \ref{genmrramg} and \ref{thmseft}, which when combined give Ramanujan's transformation \eqref{rt} (or equivalently \eqref{fwttrans}) for an integral analogue of theta functions, thus giving a better understanding of Ramanujan's transformation. Also, the results in Theorem \ref{ramtran} could have been obtained directly without having to resort to Theorems \ref{thmtrn0m4} - Corollary \ref{corexact} and Theorem \ref{appseft}. However, obtaining Theorem \ref{ramtran} from these theorems is useful since they give us many interesting results which otherwise would not have been revealed. For example, one could have proved \eqref{exact0} directly through \eqref{dires}. However, proving it through \eqref{exact} and \eqref{exact1} gives those non-trivial integral evaluations in addition. 

In light of the existence of the integral on the extreme right side of \eqref{mrramg} which equals two sides of the first error function transformation, it is natural to ask if a similar such integral exists for the second error function transformation in \eqref{thmseft}. We were unable to find such an integral and so we leave it as an open problem. However, it is important to state here the difficulty in finding this integral, if it exists.

If we reverse the steps used in proving the equality of the extreme sides of the transformation \eqref{mrramg} in Theorem \ref{genmrramg}, we notice that a crucial step was to use the integral representation for the error function given in \eqref{errapp}. However, employing the same method to the left-hand side of \eqref{thmseft} does not help as the error function there does not cancel with the integral $\displaystyle\frac{2}{\pi}\int_{0}^{\infty}\frac{e^{-\pi\alpha^2x^2+2\pi x}\sin\left(\sqrt{\pi}\alpha xz\right)}{x}\, dx$. Another reason why this is a difficult problem is that while the Mellin transform of $e^{-ax^2}\sin bx$ is essentially just a ${}_1F_{1}$ (see Lemma \ref{corinvsinl}), that of $e^{-ax^2-cx}\sin bx$ involves parabolic cylinder functions \cite[p.~503, formula 3.953, no. 1]{grn}.

We now explain the significance of this integral, provided it exists. As remarked in the introduction, subtracting the first error function transformation in \eqref{mrramg} from the second given in \eqref{seft2} leads to Ramanujan's transformation \eqref{rt} for, what is called, an integral analogue of the Jacobi theta function. The corresponding transformation for the Jacobi theta function, which has an integral involving $\Xi(t)$ equal to it \cite[Theorem 1.2]{dixthet} as well, is 
\begin{align}\label{eqsym0}
&\sqrt{\alpha}\bigg(\frac{e^{-\frac{z^2}{8}}}{2\alpha}-e^{\frac{z^2}{8}}\sum_{n=1}^{\infty}e^{-\pi\alpha^2n^2}\cos(\sqrt{\pi}\alpha nz)\bigg)&=\sqrt{\beta}\bigg(\frac{e^{\frac{z^2}{8}}}{2\beta}-e^{-\frac{z^2}{8}}\sum_{n=1}^{\infty}e^{-\pi\beta^2n^2}\cos(i\sqrt{\pi}\beta nz)\bigg)\nonumber\\
&=\frac{1}{\pi}\int_{0}^{\infty}\frac{\Xi(t/2)}{1+t^2}\nabla\left(\alpha,z,\frac{1+it}{2}\right)\, dt,
\end{align}
where the $\nabla$ function is defined in \eqref{nabla}. 

The equality of the extreme left and right sides of the special case $z=0$ of the above identity was used by Hardy \cite[Eqn.(2)]{hardyfrench} to prove that infinitely many non-trivial zeros of $\zeta(s)$ lie on the critical line Re $s=\frac{1}{2}$. Thus if an integral involving $\Xi(t)$ equal to both sides of \eqref{rt} is found, then this integral analogue of Hardy's formula may be used to obtain more information on the non-trivial zeros of $\zeta(s)$. However, this requires us to first obtain an integral involving $\Xi(t)$ equal to the two sides of \eqref{seft2}. 

Further, since our results involve an extra parameter $z$, it may be important to see what else about $\zeta(s)$, or some generalization of it, could be extracted from them. It would also be worth further studying these two error function transformations from the point of view of further applications in analytic number theory.

\textbf{Remark.} Hardy \cite{ghh} conjectured that Ramanujan's formula \eqref{mrram} may also be used for proving the infinitude of the zeros of $\zeta(s)$ on the critical line. However, we believe that it is not this formula but rather the special case $z=0$ of the identity which has an integral involving $\Xi(t)$ equal to both sides of \eqref{rt} which leads to the existence of infinitely many zeros.\\

2. Consider the transformation in \eqref{mrram} and its equivalent version \eqref{morpara} given by Mordell. Let $q=e^{i\pi w}$, Im $w>0$, and let
\begin{equation}
\Lambda(w):=\sum_{n=1}^{\infty}F(n)q^n,
\end{equation}
where $F(D)$ denotes the number of classes of positive definite binary quadratic forms $ax^2+2hxy+by^2$ with $a, b$ not both even, and determinant $-D$. Then Mordell \cite[Equation (2.18)]{mordell2} proved that
\begin{align}
\int_{0}^{\infty}\frac{xe^{\pi iwx^2}}{e^{2\pi x}-1}\, dx=-\frac{i}{4\pi w}-\Lambda(w)+\frac{\sqrt{-iw}}{w^2}\Lambda\left(-\frac{1}{w}\right)+\frac{1}{8}\left(\sum_{n=-\infty}^{\infty}e^{i\pi n^2 w}\right)^3,
\end{align}
so that with $\alpha^2=-iw$, we have for Re $\alpha^2>0$,
\begin{align}
\int_{0}^{\infty}\frac{xe^{-\pi\alpha^2x^2}}{e^{2\pi x}-1}\, dx=\frac{-1}{4\pi\alpha^2}-\sum_{n=1}^{\infty}F(n)e^{-\pi n\alpha^2}-\frac{1}{\alpha^3}\sum_{n=1}^{\infty}F(n)e^{-\pi n/\alpha^2}+\frac{1}{8}\left(\sum_{n=-\infty}^{\infty}e^{-\pi\alpha^2n^2 }\right)^3.
\end{align}
It will be interesting whether the above result admits a generalization when we work with the integral in \eqref{mrramg}.\\

3. For a fixed $z\in\mathbb{C}$, consider the integral
\begin{equation}\label{genexact}
\int_{-\infty}^{\infty}\frac{xe^{-\a x^2}}{e^{2\pi x}-1}{}_1F_{1}\left(z;\frac{3}{2};2\a x^2\right)\, dx.
\end{equation}

Using the asymptotic expansion of the confluent hypergeometric function \cite[p.~508, Equation 13.5.1]{stab}, it can be seen that as $|x|\to\infty$,
\begin{equation*}
{}_1F_{1}\left(z;\frac{3}{2};2\a x^2\right)\sim\frac{\sqrt{\pi}}{2}\left(\frac{e^{i\pi z}(2\a x^2)^{-z}}{\G\left(\frac{3}{2}-z\right)}+\frac{e^{2\a x^2}(2\a x^2)^{z-\frac{3}{2}}}{\G(z)}\right).
\end{equation*}
Note that because of the presence of $e^{2\a x^2}$ in the second expression of the asymptotic expansion, and since $\a>0$, the only way the integral in \eqref{genexact} can converge is if this expression is annihilated by $\G(z)$. This happens only when $z$ is a non-positive integer. This leads us to consider two cases based on the parity of such $z$.\\ 

\textbf{Case 1: $z$ is a non-positive even integer.} Note that for $\a$ either very small or very large, the integral in \eqref{gennearg} is nicely approximated by the expression on its right side, as in this case $\b$ is respectively very large or very small. However, the case $\a=\b=\pi$ is the worst in terms of approximating this integral, i.e., the integral
\begin{equation*}
\int_{0}^{\infty}\frac{xe^{-\pi x^2}}{e^{2\pi x}-1}{}_1F_{1}\left(-2k;\frac{3}{2};2\pi x^2\right)\, dx,
\end{equation*} 
since then $\a$ and $\b$ are not only of the same order of magnitude but in fact equal. Table 1 below shows how the integral in \eqref{gennearg} is approximated by the right side of \eqref{gennearg} for some small values of  $\a$. (The calculations in this table are done for the identity obtained by dividing both sides of \eqref{gennearg} by $\a^{3/4}$. They are performed in \emph{Mathematica}.)\\

\textbf{Case 2: $z$ is a non-positive odd integer.} When $\a=\pi$ in the integral \eqref{genexact}, we have shown in \eqref{exact0} that it is equal to zero.\\

Thus there is a trade-off in that \eqref{gennearg} has no restriction on $\a$ (except $\a>0$) but is an approximation, where as we can exactly evaluate the integral \eqref{genexact}, but only for a specific value of $\a$, i.e., when $\a=\pi$. This leads us to two open questions:

\textbf{Question 1:} Find an exact evaluation of $\displaystyle\int_{0}^{\infty}\frac{xe^{-\pi x^2}}{e^{2\pi x}-1}{}_1F_{1}\left(-2k;\frac{3}{2};2\pi x^2\right)\, dx$ for $k\in\mathbb{Z^{+}}\cup\{0\}$.\\

\textbf{Question 2:} Find an exact evaluation of, or at least an approximation to, the integral $\displaystyle\int_{0}^{\infty}\frac{xe^{-\a x^2}}{e^{2\pi x}-1}{}_1F_{1}\left(-2k-1;\frac{3}{2};2\a x^2\right)\, dx$ when $\a\neq\pi$ is a positive real number and $k\in\mathbb{Z^{+}}\cup\{0\}$.\\

It is interesting to note that in Theorem \ref{thmrgenrchi}, the integral involving $\Xi(t,\chi)$ involves the $\Delta$ function when $\chi$ is even and the $\nabla$ function when $\chi$ is odd, which is exactly opposite of what happens in Theorems 1.3, 1.4 and 1.5 in \cite{drz3}. Besides the fact that, in doing so, one can explicitly evaluate the Mellin transforms and that one \emph{does} get what one is looking for, is there some intrinsic reason behind this reversal?\\
\begin{table}[ht]
\caption{ Both side of \eqref{gennearg} (after dividing throughout by $\a^{3/4}$)}
\begin{center}
\renewcommand{\arraystretch}{1.5}
{\tiny
\begin{tabular}{c|c|c|c|c|c|c|c|c|c|c}
    \hline
    \multicolumn{1}{|c|}{$\alpha$} & \multicolumn{2}{|c|}{$.0000009$}  &\multicolumn{2}{|c|}{$.000007$} &  \multicolumn{2}{|c|}{$1.5$} & \multicolumn{2}{|c|}{$2.378$}& \multicolumn{2}{|c|}{$9361.79$}\\
    \hline
 \multicolumn{1}{|c|}{k} & \multicolumn{1}{|c|}{LHS} & \multicolumn{1}{|c|}{RHS} & \multicolumn{1}{|c|}{LHS} & \multicolumn{1}{|c|}{RHS} & \multicolumn{1}{|c|}{LHS} & \multicolumn{1}{|c|}{RHS} & \multicolumn{1}{|c|}{LHS} & \multicolumn{1}{|c|}{RHS} & \multicolumn{1}{|c|}{LHS} & \multicolumn{1}{|c|}{RHS} \\
 \hline
 \multicolumn{1}{|c|}{1}   & \multicolumn{1}{|c|}{$259259$} & \multicolumn{1}{|c|}{$259259$} & \multicolumn{1}{|c|}{$33333.4$} & \multicolumn{1}{|c|}{$33333.4$} & \multicolumn{1}{|c|}{$.212975$} & \multicolumn{1}{|c|}{$.210775$}& \multicolumn{1}{|c|}{$0.1483410$} & \multicolumn{1}{|c|}{$0.1465060$} & \multicolumn{1}{|c|}{$0.00136109$} & \multicolumn{1}{|c|}{$0.001361096$}\\
 \hline
 \multicolumn{1}{|c|}{2}   & \multicolumn{1}{|c|}{$188713$} & \multicolumn{1}{|c|}{$188713$}  & \multicolumn{1}{|c|}{$24263.1$} & \multicolumn{1}{|c|}{$24263.1$} & \multicolumn{1}{|c|}{$.162014$} & \multicolumn{1}{|c|}{$.161821$}& \multicolumn{1}{|c|}{$0.112982$} & \multicolumn{1}{|c|}{$0.112883$}& \multicolumn{1}{|c|}{$0.000990862$} & \multicolumn{1}{|c|}{$0.000990862$}\\
 \hline
 \multicolumn{1}{|c|}{3}   & \multicolumn{1}{|c|}{$154475$} & \multicolumn{1}{|c|}{$154475$}  & \multicolumn{1}{|c|}{$19861.2$} & \multicolumn{1}{|c|}{$19861.2$} & \multicolumn{1}{|c|}{$.135921$} & \multicolumn{1}{|c|}{$.137363$}& \multicolumn{1}{|c|}{$0.0948065$} & \multicolumn{1}{|c|}{$0.0960151$}& \multicolumn{1}{|c|}{$0.000811187$} & \multicolumn{1}{|c|}{$0.000811187$}\\
 \hline
  \multicolumn{1}{|c|}{4}  & \multicolumn{1}{|c|}{$133517$} & \multicolumn{1}{|c|}{$133517$} & \multicolumn{1}{|c|}{$17166.6$} & \multicolumn{1}{|c|}{$17166.6$} & \multicolumn{1}{|c|}{$.11939$} & \multicolumn{1}{|c|}{$.122057$} & \multicolumn{1}{|c|}{$0.0832805$} & \multicolumn{1}{|c|}{$0.085431$} & \multicolumn{1}{|c|}{$0.000701201$} & \multicolumn{1}{|c|}{$0.000701201$} \\
 \hline
 \multicolumn{1}{|c|}{5}    & \multicolumn{1}{|c|}{$119074$} & \multicolumn{1}{|c|}{$119074$} & \multicolumn{1}{|c|}{$15309.6$} & \multicolumn{1}{|c|}{$15309.6$} & \multicolumn{1}{|c|}{$.107718$} & \multicolumn{1}{|c|}{$.111318$}& \multicolumn{1}{|c|}{$0.0751402$} & \multicolumn{1}{|c|}{$0.07799044$} & \multicolumn{1}{|c|}{$0.000625405$} & \multicolumn{1}{|c|}{$0.000625405$}\\
 \hline
 \multicolumn{1}{|c|}{6}    & \multicolumn{1}{|c|}{$108375$} & \multicolumn{1}{|c|}{$108375$}  & \multicolumn{1}{|c|}{$13934$} & \multicolumn{1}{|c|}{$13934$} & \multicolumn{1}{|c|}{$.0989131$} & \multicolumn{1}{|c|}{$.103239$}& \multicolumn{1}{|c|}{$0.0689983$} & \multicolumn{1}{|c|}{$0.0723852$}& \multicolumn{1}{|c|}{$0.000569256$} & \multicolumn{1}{|c|}{$0.000569256$}\\
 \hline
  \multicolumn{1}{|c|}{7}    & \multicolumn{1}{|c|}{$100053$} & \multicolumn{1}{|c|}{$100053$}  & \multicolumn{1}{|c|}{$12864$} & \multicolumn{1}{|c|}{$12864$} & \multicolumn{1}{|c|}{$.091965$} & \multicolumn{1}{|c|}{$.096872$}& \multicolumn{1}{|c|}{$0.0641517$} & \multicolumn{1}{|c|}{$0.0679618$}& \multicolumn{1}{|c|}{$0.000525582$} & \multicolumn{1}{|c|}{$0.000525583$}\\
 \hline
  \multicolumn{1}{|c|}{8}    & \multicolumn{1}{|c|}{$93348.4$} & \multicolumn{1}{|c|}{$93348.4$}  & \multicolumn{1}{|c|}{$12002$} & \multicolumn{1}{|c|}{$12002$} & \multicolumn{1}{|c|}{$.0863014$} & \multicolumn{1}{|c|}{$.0916811$}& \multicolumn{1}{|c|}{$0.060201$} & \multicolumn{1}{|c|}{$0.0643522$}& \multicolumn{1}{|c|}{$0.000490396$} & \multicolumn{1}{|c|}{$0.000490397$}\\
 \hline
  \multicolumn{1}{|c|}{9}    & \multicolumn{1}{|c|}{$87801.4$} & \multicolumn{1}{|c|}{$87801.4$} & \multicolumn{1}{|c|}{$11288.8$} & \multicolumn{1}{|c|}{$11288.8$} & \multicolumn{1}{|c|}{$.0815698$} & \multicolumn{1}{|c|}{$.0873407$}& \multicolumn{1}{|c|}{$0.0569004$} & \multicolumn{1}{|c|}{$0.0613316$} & \multicolumn{1}{|c|}{$0.000461286$} & \multicolumn{1}{|c|}{$0.000461287$}\\
 \hline
 \multicolumn{1}{|c|}{10}    & \multicolumn{1}{|c|}{$83116.1$} & \multicolumn{1}{|c|}{$83116.1$}  & \multicolumn{1}{|c|}{$10686.4$} & \multicolumn{1}{|c|}{$10686.4$} & \multicolumn{1}{|c|}{$.0775398$} & \multicolumn{1}{|c|}{$.0836389$} & \multicolumn{1}{|c|}{$0.0540892$} & \multicolumn{1}{|c|}{$0.0587534$}& \multicolumn{1}{|c|}{$0.000436698$} & \multicolumn{1}{|c|}{$0.000436698$}\\
 \hline
\end{tabular}}
\end{center}
\label{tab:multicol}
\end{table}
%\begin{tabular}{c|c|c|c|c|c|c|c|c|c}
% 2.33333\times 10^6 & 1.69841\times 10^6 & 1.39028\times 10^6 & 1.20165\times 10^6 & 1.07167\times 10^6 & 975375. & 900477. & 840135. & 790212. & 748044. \\
% 466667. & 339683. & 278055. & 240331. & 214334. & 195075. & 180095. & 168027. & 158043. & 149609. \\
% 259259. & 188713. & 154475. & 133517. & 119074. & 108375. & 100053. & 93348.4 & 87801.4 & 83116.1 \\
% 233333. & 169841. & 139028. & 120166. & 107167. & 97537.6 & 90047.8 & 84013.5 & 79021.3 & 74804.5 \\
% 33333.4 & 24263.1 & 19861.2 & 17166.6 & 15309.6 & 13934. & 12864. & 12002. & 11288.8 & 10686.4 \\
% 25926. & 18871.3 & 15447.6 & 13351.8 & 11907.5 & 10837.6 & 10005.4 & 9334.91 & 8780.22 & 8311.68 \\
%\end{tabular}
\begin{center}
\textbf{Acknowledgements}
\end{center}
The first author is funded in part by the grant NSF-DMS 1112656 of Professor Victor H. Moll of Tulane University and sincerely thanks him for the support.


\begin{thebibliography}{00}
\bibitem{stab}
M.~Abramowitz and I.~A.~Stegun, Handbook of Mathematical Functions, with Formulas, Graphs, and Mathematical Tables, 9th edition, Dover publications, New York, 1961.

\bibitem{andmor}
G.~E.~Andrews, \emph{Mordell integrals and Ramanujan's ``lost'' notebook}, Analytic Number Theory (Philadelphia, PA, 1980), Lecture Notes in Mathematics, vol. 899, 1981, 10--48.

\bibitem{geabcbrln4}
G.~E.~Andrews and B.~C.~Berndt, \emph{Ramanujan's Lost Notebook. Part IV}, Springer 2013. 

\bibitem{bellman}
R.~Bellman, \emph{A Brief Introduction to Theta Functions}, Holt, Rinehart and Winston, New York, 1961; reprint Dover publications, 2013.

\bibitem{bcbspfunc}
B.C.~Berndt, Periodic Bernoulli numbers, summation formulas
and applications, in: Theory and Application of Special
Functions, R.A.~Askey, ed., Academic Press, New York, 1975,
pp.~143--189.

\bibitem{bds}
B.C.~Berndt, A.~Dixit and J.~Sohn, \emph{Character analogues of theorems of Ramanujan, Koshliakov and Guinand}, Adv. Appl. Math.~\textbf{46} (2011), 54--70 (Special issue in honor of Dennis Stanton).

\bibitem{br}
B.C.~Berndt and R.J.~Evans, \emph{Some elegant approximations and asymptotic formulas of Ramanujan}, J.~Comput.~Appl.~Math. \textbf{37} (1991), no. 1-3, 35--41.

\bibitem{berndtxu}
B.~C.~Berndt and P.~Xu, \emph{An integral analogue of theta functions and Gauss sums in Ramanujan's Lost Notebook}, Math. Proc. Cambridge Philos. Soc.~\textbf{147} (2009), 257--265.

\bibitem{cr}
B.~Chern and R.~C.~Rhoades, \emph{The Mordell integral, quantum modular forms and mock Jacobi forms}, submitted for publication.

\bibitem{con}
J.B.~Conway, \emph{Functions of One Complex Variable}, 2nd ed.,
Springer, New York, 1978.

\bibitem{dav}
H.~Davenport, \emph{Multiplicative Number Theory}, 3rd ed., Grad. Texts in Math. \textbf{74}, Springer, New York, 2000.

\bibitem{dixthet}
A.~Dixit, \emph{Analogues of the general theta transformation formula}, Proc. Roy. Soc. Edinburgh, Sect. A, \textbf{143} (2013), 371--399.

\bibitem{dixitmoll}
A.~Dixit and V.~H.~Moll, \emph{Self-reciprocal functions, powers of the Riemann zeta function and modular-type transformations}, J. Number Thy.~\textbf{147} (2015), 211--249.

\bibitem{drz3}
A.~Dixit, A.~Roy and A.~Zaharescu, \emph{Riesz-type criteria and theta transformation analogues}, submitted for publication.

\bibitem{tti}
A.~Erd\'elyi, W.~Magnus, F.~Oberhettinger and F.G.~Tricomi, \emph{Tables of Integral Transforms}, Vol. 1, McGraw-Hill, New York, 1954.

\bibitem{glashier1} 
J.~W.~L.~Glashier, \emph{XXXII. On a class of definite integrals}, The London, Edinburgh, and Dublin Philosophical Magazine and Journal of Science,~\textbf{42}, No. 280 (1871), 294--302.

\bibitem{glashier2}
J.~W.~L.~Glashier, \emph{LIV. On a class of definite integrals-Part II}, The London, Edinburgh, and Dublin Philosophical Magazine and Journal of Science,~\textbf{42}, No. 282 (1871), 421--436.

\bibitem{grn}
I.~S.~Gradshteyn and I.~M.~Ryzhik, eds., \emph{Table of Integrals,
Series, and Products}, 7th ed., Academic Press, San Diego, 2007.

\bibitem{hardyfrench}
G.~H.~Hardy, \emph{Sur les z\'{e}ros de la fonction $\zeta(s)$ de Riemann}, CR Acad. Sci. Paris~\textbf{158} (1914), 1012--1014.

\bibitem{ghh}
G.H.~Hardy, \emph{Note by G.H.~Hardy on the preceding paper}, Quart.~J.~Math.~\textbf{46} (1915),
260--261.

\bibitem{hiary1}
G.~A.~Hiary, \emph{Fast methods to compute the Riemann zeta function}, Ann. Math.~\textbf{174} (2011), 891--946.

\bibitem{hiary2}
G.~A.~Hiary, \emph{A nearly-optimal method to compute the truncated theta function, its derivatives,
and integrals}, Ann. Math.~\textbf{174} (2011), 859--889.

\bibitem{jl}
E.~Jahnke, F.~Emde and F.~L\"osch, \emph{Tables of higher functions}, 6th ed., McGraw-Hill, New York, 1960.

\bibitem{koshli}
N.S.~Koshliakov, \emph{Note on some infinite integrals}, C. R. (Dokl.) Acad. Sci. URSS ~\textbf{2} (1936), 247--250.

\bibitem{kosh1937}
N.S.~Koshliakov, \emph{On a transformation of definite integrals
and its application to the theory of Riemann's function
$\zeta(s)$}, Comp.~Rend.~(Doklady) Acad.~Sci.~URSS \textbf{15}
(1937), 3--8.

%\bibitem{kramp}
%C.~Kramp, \emph{Analyse des R\'{e}fractions Astronomiques et Terrestres}, Strasbourg, 1798.

\bibitem{kronecker1}
L.~Kronecker, \emph{Summirung der Gausschen Reihen $\sum_{h=0}^{h=n-1}e^{\frac{2h^2\pi i}{n}}$}, J. Reine Angew. Math.,~\textbf{105} (1889), 267--268.

\bibitem{kronecker2}
L.~Kronecker, \emph{Bemerkungen uber die Darstellung von Reihen durch integrale. (Fortsetzung.)}, J. Reine Angew. Math.,~\textbf{105} (1889), 345--354.

\bibitem{kuznetsov}
A.~Kuznetsov, \emph{Computing the truncated theta function via Mordell integral}, to appear in \emph{Mathematics of Computation}, \url{http://dx.doi.org/10.1090/mcom/2953}.

\bibitem{lerch1}
M.~Lerch, \emph{Bemerkungen zur Theorie der elliptischen Funktionen}, Rozpravyceske Akademie cisare Frantiska
Josefa pro vedy slovesnost, a umeni (II Cl) Prag. (Bohmisch) I Nr. (24) 1892. Abstract in the
Jahrbuch \"{u}ber die Fortschritte der Mathematik, 24 (1892), 442--445

\bibitem{lerch2}
M.~Lerch, \emph{Beitr\"{a}ge zur Theorie der elliptischen Funktionen, unendlichen Reihen und bestimmter Integrale}, Rozpravyceske Akademie cisare Frantiska
Josefa pro vedy slovesnost, a umeni (II Cl) Prag. (Bohmisch) II Nr. 23 (1983).

\bibitem{lerch3}
M.~Lerch, \emph{Zur Theorie der Gausschen Summen} Math. Ann.,~\textbf{57} (1903), 554--567.

\bibitem{mordell1}
L.~J.~Mordell, \emph{The value of the definite integral $\displaystyle\int_{-\infty}^{\infty}\frac{e^{at^2+bt}}{e^{ct}+d}\, dt$}, Q.~J.~Math.~\textbf{68} (1920), 329--342.

\bibitem{mordell2}
L.~J.~Mordell, \emph{The definite integral $\displaystyle\int_{-\infty}^{\infty}\frac{e^{at^2+bt}}{e^{ct}+d}\, dt$ and the analytic theory of numbers}, Acta Math.~\textbf{61} (1933), 323--360.

\bibitem{ober}
F.~Oberhettinger, \emph{Tables of Mellin Transforms}, Springer-Verlag, New York, 1974.

\bibitem{oloa}
O.~Oloa, \emph{On a series of Ramanujan}, in
{\it Gems in Experimental Mathematics}, T.~Amdeberhan and
V.~H.~Moll, eds., Contemp.~Math., Vol.~517, American Mathematical
Society, Providence, RI, 2010, pp.~305--311.

\bibitem{olver}
F.W.J.~Olver, \emph{Asymptotics and Special Functions}, Academic Press, New York, 1974.

\bibitem{nist}
F.~W.~J.~Olver, D.~W.~Lozier, R.~F.~Boisvert, and C.~W.~Clark, editors. NIST Handbook of
Mathematical Functions. Cambridge University Press, 2010.

\bibitem{kp}
R.B.~Paris and D.~Kaminski, \emph{Asymptotics and Mellin-Barnes Integrals},  Encyclopedia of Mathematics and its Applications, 85. Cambridge University Press, Cambridge, 2001.

\bibitem{rain}
E.D.~Rainville, \emph{Special Functions}, The Macmillan company, New York, 1960.

\bibitem{sdi}
S.~Ramanujan, \emph{Some definite integrals}, Mess. Math.~\textbf{44}, (1915), 10--18.

\bibitem{riemann}
S.~Ramanujan, \emph{New expressions for Riemann's functions $\xi(s)$ and $\Xi(t)$}, Quart. J. Math.~\textbf{46} (1915), 253--261.

\bibitem{ram1915}
S.~Ramanujan, \emph{Some definite integrals connected with Gauss's sums}, Mess. Math.~\textbf{44} (1915), 75--85.

\bibitem{ram1919}
S.~Ramanujan, \emph{Some definite integrals}, J. Indian Math. Soc.~\textbf{11} (1919), 81--87.

\bibitem{ramcol}
S.~Ramanujan, Collected Papers, Cambridge University Press, Cambridge, 1927, reprinted by Chelsea, New York, 1962, reprinted by the Amer. Math.
Soc., Providence, RI, 2000.

\bibitem{ramnot}
S.~Ramanujan, \emph{Notebooks}, Tata Institute of Fundamental Research, Bombay, 1957 (2 volumes); reprint 2012.

\bibitem{lnb}
S.~Ramanujan, \emph{The Lost Notebook and Other Unpublished
Papers}, Narosa, New Delhi, 1988.

\bibitem{siegel}
C.~L.~Siegel, \emph{\"{U}ber Riemanns Nachlass zur analytischen Zahlentheorie.} Quellen und Studien zu Geschichte de Mathematik, Astronomie und Physik,~\textbf{2} (1933), 45--80.

\bibitem{temme}
N.M.~Temme, \emph{Special functions: An introduction to the classical functions of mathematical physics}, Wiley-Interscience Publication, New York, 1996.

\bibitem{titch}
E.~C.~Titchmarsh, \emph{The Theory of the Riemann Zeta Function}, Clarendon
Press, Oxford, 1986.

\bibitem{zwegers}
S.~P.~Zwegers, \emph{Mock Theta Functions}, PhD thesis, Universiteit Utrecht, 2001.\\
\url{http://arxiv.org/pdf/0807.4834v1.pdf} 
\end{thebibliography}
\end{document}